\documentclass[11pt]{amsart}

\usepackage{amsmath,amsthm, amscd, amssymb, amsfonts}
\usepackage{amsthm}
\usepackage[all]{xy}

\usepackage[ansinew]{inputenc}
\usepackage{graphicx,fancyhdr}

\usepackage[dvips, dvipsnames, usenames]{color}
\newtheorem{theorem}{Theorem}[section]

\newtheorem{lemma}[theorem]{Lemma}
\newtheorem{coro}[theorem]{Corollary}

\newtheorem{prop}[theorem]{Proposition}

\theoremstyle{definition}
\newtheorem{definition}[theorem]{Definition}
\newtheorem{example}[theorem]{Example}

\theoremstyle{remark}
\newtheorem{remark}[theorem]{Remark}

\newcommand{\vi}{\textbf{(i)} }
\newcommand{\vii}{\textbf{(ii)} }
\newcommand{\viii}{\textbf{(iii)} }
\newcommand{\viv}{\textbf{(iv)} }

\newcommand{\bq}{\mathbf{q} }

\newcommand{\ydg}{^{\ku\Gamma}_{\ku\Gamma}\mathcal{YD}}
\newcommand{\ydl}{^{\ku\Lambda}_{\ku\Lambda}\mathcal{YD}}

\newcommand\Z{\mathbb Z}
\newcommand\N{\mathbb N}

\def\cI{\mathcal{I}}

\def\cB{\mathcal{B}}
\def\cU{\mathcal{U}}

\def\cX{\mathcal{X}}
\def\cG{\mathcal{G}}
\def\cC{\mathcal{C}}
\def\cW{\mathcal{W}}
\def\cR{\mathcal{R}}
\def\cS{\mathcal{S}}
\def\cD{\mathcal{D}}
\def\cK{\mathcal{K}}
\def\cH{\mathcal{H}}
\newcommand\ad{\operatorname{ad}}
\newcommand{\ku}{ \mathbf{k}}
\newcommand{\eps}{\varepsilon}

\newcommand\cop{\operatorname{cop}}
\newcommand\Hom{\operatorname{Hom}}
\newcommand\id{\operatorname{id}}
\def\ot{\otimes}
\def\u{\mathfrak{u}}
\def\Zt{\Z^{\theta}}

\newcommand{\Eb}{\underline E}
\newcommand{\Fb}{\underline F}
\newcommand{\Kb}{\underline K}
\newcommand{\Lb}{\underline L}
\newcommand{\qb}{\underline q}

\newcommand{\schi}{{s_p^*\chi}}
\newcommand{\de}{\Delta}
\newcommand{\En}{\mathbf{E}}

\def\pf{\begin{proof}}
\def\epf{\end{proof}}

  \makeatletter
  \@addtoreset{equation}{section}
  \makeatother

\newcommand{\trDeltatwo}{{\trDelta^{(2)}}}

\newcommand{\trrho}{\rho}
\newcommand{\trV}{V}
\newcommand{\truchi}{{\u(\chi)}}
\newcommand{\rmEnd}{{\rm{End}}}
\newcommand{\trv}{v}
\newcommand{\trLambda}{\Lambda}
\newcommand{\trupchi}{{\u^+(\chi)}}
\newcommand{\trumchi}{{\u^-(\chi)}}
\newcommand{\truzerochi}{\u^0(\chi)}

\newcommand{\trugeqchi}{{\u^{\geq 0}(\chi)}}
\newcommand{\truleqchi}{{\u^{\leq 0}(\chi)}}

\newcommand{\trK}{K}
\newcommand{\trL}{L}
\newcommand{\trE}{E}
\newcommand{\trF}{F}

\newcommand{\trtheta}{\theta}

\newcommand{\mcF}{{\mathcal{F}}}
\newcommand{\rmid}{{\mathrm{id}}}

\newcommand{\mcC}{{\mathcal{C}}}

\newcommand{\Nzerotheta}{{\N_0^\theta}}

\newcommand{\treta}{\eta}
\newcommand{\trf}{f}
\newcommand{\rmGL}{{\mathrm{GL}}}
\newcommand{\trchi}{\chi}
\newcommand{\tral}{\alpha}
\newcommand{\trbeta}{\beta}
\newcommand{\trgamma}{\gamma}

\newcommand{\trR}{R}

\newcommand{\trS}{S}
\newcommand{\trDelta}{\Delta}

\newcommand{\trtau}{\tau}

\newcommand{\hattretageq}{{{\hat \treta}^\geq}}
\newcommand{\hattretaleq}{{{\hat \treta}^\leq}}

\newcommand{\cheX}{{\check{X}}}

\begin{document}

\title[$R$-matrix of quantum doubles of Nichols algebras]{The $R$-matrix of quantum doubles of Nichols algebras of diagonal type}
\author[Angiono, Yamane]{Iv\'an Angiono}
\address{FaMAF-CIEM (CONICET), Universidad Nacional de C\'ordoba,
Medina A\-llen\-de s/n, Ciudad Universitaria (5000) C\' ordoba, Rep\'
ublica Argentina.}
\email{angiono@famaf.unc.edu.ar}
\author[]{Hiroyuki Yamane}
\address{University of Toyama,
Faculty of Science,
Gofuku 3190, Toyama-shi, Toyama 930-8555, Japan}
\email{hiroyuki@sci.u-toyama.ac.jp}

\thanks{The work of I. A. was partially supported by
CONICET, FONCyT-ANPCyT, Secyt (UNC).}

\begin{abstract}
Let $H$ be the quantum double of a Nichols algebra of diagonal type.
We compute the $R$-matrix of 3-uples of modules for general finite-dimensional
highest weight modules over $H$.

We calculate also a multiplicative formula for the universal $R$-matrix when $H$ is finite dimensional.
\end{abstract}

\maketitle

\section{Introduction}

A remarkable property of quantum groups, introduced by Drinfeld and Jimbo in the eighties, is the existence of an $R$- matrix
for their categories of modules. This $R$-matrix is related with the existence of solutions of the Yang-Baxter equation. An
explicit formula for the universal $R$-matrix of quantum groups was obtained in the nineties \cite{KR,LS,R-PBWThm}, and
extended to quantized enveloping superalgebras \cite{KhT,Y} of finite-dimensional Lie superalgebras.

We can deduce the existence of this $R$-matrix for quantized enveloping (super)algebras because they can be obtained as quotients
of quantum doubles of bosonizations of the positive part by group algebras, and these quantum doubles are quasi-triangular.

A natural generalization of the positive part of these quantized enveloping (super)algebras are the Nichols algebras of diagonal type
\cite{ASpointedHA}. They admit a root system and a Weyl groupoid \cite{HY-weylgpd, HS} controlling the structure of these algebras.
Moreover, the classification of these Nichols algebras with finite root system includes (properly) the positive part of quantized
enveloping algebras of finite dimensional contragradient Lie superalgebras and simple Lie algebras. It is natural then to ask for a formula of the $R$-matrix
in this general context. We answer this question for the subfamily of finite-dimensional representations with a highest weight in a
general context, and obtain an explicit formula for the universal $R$-matrix when the Nichols algebra is finite-dimensional.

\medskip

Although Nichols algebras appeared as an important tool for the classification of finite dimensional pointed Hopf algebras \cite{ASpointedHA}, they have become very attractive for another fields of mathematics. In particular they are related with conformal field theories. Indeed they give place to logarithmic examples \cite{ST - screening,ST - log CFT,ST - repr Uq sl}. Starting from non-semisimple (logarithmic) CFT and the screening operators, we can obtain a braided Hopf algebra which is a Nichols algebra \cite{ST - screening}. Then it becomes interesting how to make a reverse construction in order to obtain new examples of vertex operators algebras and the corresponding conformal field theories. This was the motivation to their study in Mathematical Physics \cite{ST - log CFT}: these authors start the translation of some elements from the Nichols algebra context to the corresponding ones needed to describe the attached vertex operator algebra. They study the category of Yetter-Drinfeld modules over the Nichols algebra into a braided category, which is exactly the category of representations of the quantum double of the bosonization of this Nichols algebra by the group algebra of a finite abelian group. They complete the computation for a particular example \cite{ST - repr Uq sl}, describing the projective modules, and they give the $R$-matrix following the present work. The $R$-matrix encodes the $M$-matrix for the dual algebra of the corresponding quantum double, which is responsible of the monodromy on the CFT language \cite{FGST}.

\medskip

The organization of the paper is as follows. In Section \ref{section:Lusztig isomorphisms} we recall definitions
and results needed for our work. They are related with quantum doubles and properties of Nichols algebras of diagonal
type. We stress the importance of the Weyl groupoid and the generalized version of root systems. In Section \ref{section:Rfnrep}
we work over arbitrary Nichols algebras of diagonal type and compute the $R$-matrix of 3-uples of finite-dimensional modules,
generalizing the results in \cite{Tanisaki}. We restrict our attention to highest weight modules, which give maybe the most
important subfamily of representations. Finally in Section \ref{section:universal Rmatrix} we compute the universal $R$-matrix
for quantum doubles of finite-dimensional Nichols algebras. The formula involves the multiplication of quantum exponentials
of root vector powers, generalizing the classical ones for quantum groups.

\subsection*{Notation} We denote by $\N$ the set of natural numbers,
and by $\N_0$ the set of non-negative integers.

Let $\ku$ be an algebraically closed field of characteristic zero. All the vector spaces, algebras and tensor products are over $\ku$.
We shall use the usual notation for $q$-combinatorial numbers: for each $q\in\ku^\times$, $n\in\N$, $0\leq k\leq n$,
\begin{align*}
(n)_q&=1+q+\ldots+q^{n-1}, & (n)_q!&=(1)_q(2)_q\cdots(n)_q,
\\ \binom{n}{k}_q&=\frac{(n)_q!}{(k)_q!(n-k)_q!}.
\end{align*}

Let $A$ be an associative algebra. Given an element $a\in A$ such
that $a^N=0$, we define the $q$-\emph{exponential}, for each $q$ which is not a root of unity, or it is a root of unity
of order $\geq N$:
\begin{equation}
\exp_q(a)=\sum_{i=0}^{N-1} \frac{a^i}{(i)_q!}.
\end{equation}

\medskip

Let $\theta\in\N$. $\{\alpha_i\}_{1\leq i\leq\theta}$ will denote the canonical $\Z$-basis of $\Z^\theta$.
Given an $\Z$-linear automorphism $s:\Z^\theta \to
\Z^\theta$ and a bicharacter $\chi : \Zt \times \Zt \to
\ku^{\times}$, $s^*\chi$ will denote the bicharacter
\begin{equation}\label{eqn:accion w sobre chi}
(s^*\chi)(\alpha,\beta):= \chi\left(s^{-1}(\alpha), s^{-1}(\beta)
\right), \qquad \alpha, \beta \in \Z^\theta.
\end{equation}

\medskip

Given a Hopf algebra $H$ with coproduct $\Delta$, we will use the classical Sweedler notation $\Delta(h)=h_1\ot h_2$, $h\in H$, and denote
$$\trDeltatwo:= (\trDelta\otimes\rmid)\circ\Delta =(\rmid\otimes\trDelta)\circ\Delta.$$
Given $R=\sum_i a_i\otimes b_i\in H\otimes H$, we set the following elements of $H\otimes H\otimes H$:
\begin{align*}
R^{(1,2)}&=\sum_i a_i\otimes b_i\otimes 1, & R^{(1,3)}&=\sum_i a_i\otimes1\otimes b_i, & R^{(2,3)}&=\sum_i 1\otimes a_i\otimes b_i.
\end{align*}

Recall that a (right) coideal subalgebra of $H$ is a subalgebra $A$ of $H$ such that $\Delta(A)\subseteq H\otimes A$.

\section[Preliminaries]{Preliminaries}\label{section:Lusztig isomorphisms}

We recall some definitions and results which will be useful in the rest of this work. They are mainly related with
quantum doubles of Hopf algebras and Nichols algebras of diagonal type.

\subsection{Skew-Hopf pairings and $R$-matrices}\label{subsection:skew-hopf pairing}

Let $A,B$ be two Hopf algebras. A \emph{skew Hopf pairing}
between $A$ and $B$ (see \cite[Section 3.2.1]{J}, \cite[Section
8.2]{KS}) is a linear map $\eta:A\otimes B\to \ku$ such that
\begin{align*}
& \eta(xx',y)=\eta(x',y_1)(x,y_2), & &\eta(x,1)=\varepsilon(x), \\
& \eta(x,yy')=\eta(x_1,y)(x_2,y'), & &\eta(1,y)=\varepsilon(y), \\
& \eta(\cS(x),y)=\eta(x,\cS^{-1}(y)), & &
\end{align*}
for all $x,x'\in A$, $y,y'\in B$. In such case, $A\otimes B$ admits
a unique structure of Hopf algebra, denoted by $\cD(A,B,\eta)$ and
called the \emph{quantum double} associated to $\eta$, such
that the morphisms $A\to A\otimes B$, $a\mapsto a\otimes 1$, $B\to
A\otimes B$, $b\mapsto 1\otimes b$ are Hopf algebra morphisms and
$$ (a\otimes 1)(1\otimes b)=a\otimes b, \qquad
(1\otimes b)(a\otimes 1)= \eta\left(a_1,\cS(b_1)\right)(a_2\otimes
b_2)\eta(a_3,b_3).$$

When $A$ is finite-dimensional and $\eta$ is not degenerate, $B$ is
identified with the Hopf algebra $A^*$. $\cD(A,B,\eta)=\cD(A)$ is
the \emph{Drinfeld double} of $A$, which admits an $R$-matrix:
\begin{equation}\label{eqn:R-matrix double}
    \cR:= \sum_{i\in I} (1\otimes b_i)\otimes (a_i\otimes 1),
\end{equation}
where $\{a_i\}_{i\in I}$, $\{b_i\}_{i\in I}$ are dual bases of $A$,
$B$: $\eta(a_i,b_j)=\delta_{ij}$.

\subsection{Weyl groupoids and convex orders on finite root systems}
We recall the definitions of Weyl groupoid and generalized root system following \cite{CH}.
Fix a set $\cX\neq \emptyset$ and a finite set $I$. Set also for each $i \in I$ a map
$r_i:\cX\rightarrow\cX$, and for each $X\in\cX$ a generalized Cartan matrix $A^X=(a^X_{ij})_{i,j\in I}$
in the sense of \cite{K}.

\begin{definition}{\cite{CH,HY-weylgpd}}
The quadruple $\cC:= \cC(I, \cX, (r_i)_{i \in I}, (A^X)_{X \in \cC})$ is a \emph{Cartan scheme} if $r_i^2=\id$
for all $i \in I$, and $a^X_{ij}=a^{r_i(X)}_{ij}$ for all $X \in \cX$ and $i,j \in I$. For each $i \in I$ and
each $X \in \cX$ set $s_i^X$ as the $\Z$-linear automorphism of $\Z^I$ determined by
\begin{align}\label{eq:def s_i}
s_i^X(\alpha_j)&=\alpha_j-a_{ij}^X\alpha_i, & j&\in I.
\end{align}
The \emph{Weyl groupoid} of $\cC$ is the groupoid $\cW(\cC)$ whose set of objects is $\cX$ and whose morphisms are generated by $s_i^X$; here we consider $s_i^X$ as an element in $\Hom(X, r_i(X))$, $i \in I$, $X \in \cX$.

Given a Cartan scheme $\cC$, and for each $X \in \cX$ a set $\de^X \subset \Z^I$, we say that $\cR:= \cR(\cC, (\de^X)_{X \in \cX} )$ is a \emph{root system of type} $\cC$ if
\begin{itemize}
  \item for all $X \in \cX$, $\de^X= (\de^X \cap \N_0^I) \cup  -(\de^X \cap \N_0^I)$. We call $\de^X_+:= \de^X \subset \N_0^I$ the set of \emph{positive roots}, and $ \de^X_-:= - \de^X_+$ the set of \emph{negative roots}.
  \item for all $i \in I$ and all $X \in \cX$, $\de^X \cap \Z \alpha_i= \{\pm \alpha_i \}$.
  \item for all $i \in I$ and all $X \in \cX$, $s_i^X(\de^X)=\de^{r_i(X)}$.
  \item set $m_{ij}^X:= |\de^X \cap (\N_0\alpha_i+\N_0 \alpha_j)|$; then $(r_ir_j)^{m_{ij}^X}(X)=(X)$ for all $i \neq j \in I$ and all $X\in \cX$.
\end{itemize}
\end{definition}
We assume that $\cW(\cC)$ is a connected groupoid: $\Hom(Y,X) \neq \emptyset$, for all $X,Y \in \cX$. For any $X \in \cX$, let
$\Hom(\cW,X):= \cup_{Y \in \cX} \Hom(Y,X)$, and $\de^{X \ re}:= \{ w(\alpha_i): \ i \in I, \ w \in \Hom(\cW,X) \}$, the set of \emph{real roots} of $X$. Clearly, $\de^{X \ re}\subseteq \de^X$, and $w(\de^X)= \de^Y$ for any $w \in \Hom(Y,X)$. We say that $\cR$ is \emph{finite} if $\de^X$ is finite for some $X\in \cX$.

Note that each $w \in \Hom(\cW,X_1)$ is a product $s_{i_1}^{X_1}s_{i_2}^{X_2} \cdots s_{i_m}^{X_m}$, where $X_j=r_{i_{j-1}} \cdots r_{i_1}(X_1)$, $i \geq 2$; we fix the notation $w= \id_{X_1} s_{i_1} \cdots s_{i_m}$ to mean that $w \in \Hom(\cW,X_1)$, because the elements $X_j \in \cX$ are determined by the previous condition. The \emph{length} of $w$ is
$$ \ell(w)= \min \{ n \in \N_0: \ \exists i_1, \ldots, i_n \in I \mbox{ such that }w=\id_X s_{i_1} \cdots s_{i_n} \}. $$

\begin{prop}{\cite[Prop. 2.12]{CH}}  \label{Prop:maxlengthCH}
Let $w=\id_X s_{i_1} \cdots s_{i_m}$, $\ell(w)=m$. The roots $\beta_j=s_{i_1} \cdots s_{i_{j-1}}(\alpha_{i_j}) \in \de^X$ are positive and pairwise different.

Moreover, if $\cR$ is finite and $w$ is an element of maximal length, then $\{ \beta_j \}= \Delta^X_+$, so all the roots are real. \qed
\end{prop}

For the last part of this subsection, assume that $\cR$ is finite.

\begin{definition}{\cite{A2}}
Given a root system $\cR$ and a fixed total order $<$ on $\Delta^X_+$, we say that it is
\emph{convex} if for each $\alpha, \beta \in \Delta^X_+$ such that $\alpha <  \beta$ and $\alpha+\beta \in \Delta^X_+$, then
$ \alpha < \alpha+\beta <\beta$. It is said \emph{strongly convex} if for each ordered subset
$\alpha_1 \leq  \ldots \leq \alpha_k$ of elements of $\Delta^X_+$ such that
$\alpha := \sum \alpha_i \in \Delta^X_+$, it holds that $\alpha_1 < \alpha <
\alpha_k$.
\end{definition}

\begin{theorem}{\cite{A2}}\label{Theorem:equivalenciasconvexo}
Given an order on $\Delta^X_+$, the following are equivalent:
\begin{enumerate}
  \item the order is convex,
  \item the order is strongly convex,
  \item the order is associated with a reduced expression of the longest element. \qed
\end{enumerate}
\end{theorem}

\subsection{Weyl groupoid of a Nichols algebra of diagonal type}\label{subsection:weyl groupoid}

For each bicharacter $\chi:\Zt\times\Zt\to\ku^{\times}$, set $q_{ij}(\chi)=\chi(\alpha_i, \alpha_j)$.
Given $1\leq i\leq\theta$, we say that $\chi$ is \emph{$i$-finite} if for all $1\leq j\neq i\leq \theta$
there exists $m\in\N_0$ such that $(m+1)_{q_{ii}}(1-q_{ii}^2q_{ij}q_{ji})=0$. In such case, define
\begin{align*}
a_{ii}^\chi&=2, & a_{ij}^\chi &=-\min\{ m\in\N_0| (m+1)_{q_{ii}}(1-q_{ii}^2q_{ji}q_{ij})=0\},
\end{align*}
and set $s_i^\chi$ as the $\Z$-linear automorphism of $\Z^\theta$ given by \eqref{eq:def s_i}.
If $\chi$ is $i$-finite for all $i$, $A^\chi=(a_{ij}^\chi)_{1\leq i,j\leq\theta}$ is the \emph{generalized Cartan matrix} associated to $\chi$.

Let $\cX$ be the set of all the bicharacters of $\Zt$. We define $r_i:\cX\to\cX$ by $r_i(\chi)=(s_i^\chi)^*\chi$ if $\chi$ is $i$-finite, or $r_i(\chi)=\chi$ otherwise. Such $r_i$'s are involutions and $\cG(\chi)$ will denote the orbit of $\chi$ by the action of the group of bijections generated by the $r_i$'s.

Note that $\cC(\chi)=\cC\left(\{1,\ldots,\theta\},\cG(\chi),(r_i)_{1\leq i\leq\theta},(C^{\upsilon})_{\upsilon\in\cG(\chi)} \right)$ is a connected Cartan scheme, see \cite{HY-weylgpd,HY-shapov}. Therefore the associated Weyl groupoid $\cW(\chi)$ is called the \emph{Weyl groupoid of $\chi$}.

There exists a close relation between the root system and the set $\cK(V)$ of graded coideal subalgebras, as it is stated in \cite{HS}. We refer the reader to this reference for the definitions of the coideal subalgebra $B(w)$ and the Duflo order.

\begin{theorem}\cite{HS}\label{Thm:HSadaptado}
For each $w \in \Hom(\cW,V)$ there exists a unique right coideal
subalgebra $B(w) \in \cK(V)$ such that its Hilbert series is
\begin{equation}\label{eqn:HilbertseriesF}
    \cH_{B(w)} = \prod_{\beta \in \Lambda^V_+(w)} \bq_{N_\beta}(X^{\beta}).
\end{equation}
Moreover, the correspondence $w \mapsto B(w)$ gives an order
preserving and order reflecting bijection between $\Hom(\cW,V)$ and
$\cK(V)$, where we consider the Duflo order over $\Hom(\cW,V)$ and
the inclusion order over $\cK(V)$; i.e. $ w_1 \leq_D w_2$ if and
only if $B(w_1) \subset B(w_2)$.\qed
\end{theorem}

\subsection{Lusztig Isomorphisms of Nichols algebras}

Set $\chi$, $(q_{ij})$ as in Subsection \ref{subsection:weyl
groupoid}. $\cU(\chi)$ will denote the algebra presented by generators $E_i$,
$F_i$, $K_i$, $K_i^{-1}$, $L_i$, $L_i^{-1}$, $1 \leq i \leq \theta$,
and relations
\begin{align*}
XY&=YX, & X,Y \in & \{ K_i^{\pm1}, L_i^{\pm1}: 1 \leq i \leq
\theta \},
\\ K_iK_i^{-1}&=L_iL_i^{-1}=1, & E_iF_j-F_jE_i&=\delta_{i,j}(K_i-L_i)
\\ K_iE_jK_i^{-1}&=q_{ij}E_j, \quad & L_iE_jL_i^{-1}&=q_{ji}^{-1}E_j,
\\ L_iF_jL_i^{-1}&=q_{ji}F_j, \quad & K_iF_jK_i^{-1}&=q_{ij}^{-1}F_j.
\end{align*}
$\cU^{+0}(\chi)$ (respectively, $\cU^{-0}(\chi)$) will denote the subalgebra
generated by $K_i$, $K_i^{-1}$ (respectively, $L_i$, $L_i^{-1}$), $1
\leq i \leq \theta$, and $\cU^0(\chi)$ will denote the subalgebra
generated by $K_i$, $K_i^{-1}$, $L_i$ and $L_i^{-1}$.
Also, $\cU^+(\chi)$ (respectively, $\cU^-(\chi)$) will denote the subalgebra
generated by $E_i$ (respectively, $F_i$), $1 \leq i \leq \theta$.

$\cU(\chi)$ is a $\Zt$-graded Hopf algebra, with graduation
determined by the following conditions:
\begin{align*}
\deg(K_i)&=\deg(L_i)=0, & \deg(E_i)&=\alpha_i, & \deg(F_i)&=-\alpha_i.
\end{align*}
$\cU(\chi)$ admits a Hopf algebra structure, with comultiplication determined by
\begin{align*}
\Delta(K_i)&=K_i \ot K_i, & \Delta(E_i)&=E_i \ot 1 + K_i \ot E_i,
\\ \Delta(L_i)&=L_i \ot L_i, & \Delta(F_i)&=F_i \ot L_i + 1 \ot F_i,
\end{align*}
and then $\eps(K_i)=\eps(L_i)=1$, $\eps(E_i)=\eps(F_i)=0$.

Note that $\cU^0(\chi)$ is isomorphic to $\ku \Z^{2\theta}$ as Hopf algebras, and
the subalgebra $\cU^{\geq0}(\chi)$ (respectively, $\cU^{\leq0}(\chi)$) generated by
$\cU^+(\chi)$, $K_i^{\pm1}$, $1 \leq i \leq \theta$, (respectively, $\cU^-(\chi)$, $L_i^{\pm1}$) is isomorphic to
$T(V)\#\ku\Zt$ (respectively, $T(V^*)\#\ku\Zt$). $\cU(\chi)$ is the associated quantum double.

Here, $\cU^+(\chi)$ is isomorphic to $T(V)$ as
braided graded Hopf algebras in the category of Yetter-Drinfeld
modules over $\ku \Z^\theta$, with actions and
coactions given by:
$$ K_i \cdot E_j=q_{ij} E_j, \qquad \delta(E_i)= K_i \ot E_i, $$
and similar equations for $F_j$, $L_i$.
$\underline \Delta$ will denote the braided comultiplication of
$\cU^+(\chi)$. As it is $\N_0$-graded, we will consider
$\underline\Delta_{n-k,k}(E)$, the component of $\underline
\Delta(E)$ in $\cU^+(\chi)_{n-k}\ot\cU^+(\chi)_k$, for each
$E\in\cU^+(\chi)$ homogeneous of degree $n$, and
$k\in\{0,1,\ldots,n\}$.

By \cite[Prop. 4.14]{H-isom}, the multiplication $ m: \cU^+(\chi)
\ot \cU^0(\chi) \ot \cU^-(\chi) \to \cU(\chi)$ is an isomorphism of
$\Z^\theta$-graded vector spaces.

\medskip

We consider some isomorphisms involving $\cU(\chi)$ \cite[Section
4.1]{H-isom}.
\begin{enumerate}
  \item[(a)] Let $\underline{a}=(a_1,\ldots,a_\theta)\in(\ku^\times)^\theta$.
There exists a unique algebra automorphism $\varphi_{\underline a}$
of $\cU(\chi)$ such that
\begin{equation}\label{eqn:varphi a}
\varphi_{\underline a}(K_i)=K_i, \quad \varphi_{\underline
a}(L_i)=L_i, \quad \varphi_{\underline a}(E_i)=a_iE_i, \quad
\varphi_{\underline a}(F_i)=a_i^{-1}F_i.
\end{equation}
  \item[(b)] There exists a unique algebra antiautomorphism $\Omega$ of $\cU(\chi)$ such
  that
\begin{equation}\label{eqn:omega}
 \Omega(K_i)=K_i, \quad \Omega(L_i)=L_i, \quad \Omega(E_i)=F_i, \quad \Omega(F_i)=E_i.
\end{equation}
It satisfies the relation $\Omega^2=\id$.
\end{enumerate}

As in \cite{HY-shapov,H-isom}, $\cI^+(\chi)$ will denote the ideal
of $\cU^+(\chi)$ such that the quotient $\cU^+(\chi)/\cI^+(\chi)$ is
isomorphic to the Nichols algebra of $V$; that is, the greatest
braided Hopf ideal of $\cU^+(\chi)$ generated by elements of degree
$\geq2$. Set $\cI^-(\chi)= \Omega(I^+(\chi))$, where $\phi_4$ is the
anti-automorphism of algebras determined by \eqref{eqn:omega}, and
$$ \u^{\pm}(\chi):= \cU^{\pm}(\chi) / \cI^{\pm}(\chi), \qquad \u(\chi):= \cU(\chi) / (\cI^-(\chi)+\cI^+(\chi)),$$
and $\u^{\geq0}(\chi)$, $\u^{\leq0}(\chi)$ the corresponding images on the quotient.
Note that $\u(\chi)$ is the quantum double of
$\u^+(\chi)\#\ku\Z^\theta$. The following result follows by
\cite[Lemma 6.5, Theorem 6.12]{H-isom}.

\begin{prop}\label{prop:non deg pairing HY}
\cite[Proposition 3.5]{HY-shapov}, \cite[Theorem 5.8]{H-isom}
There exists a unique non-degenerate skew-Hopf pairing $\eta:\u^+(\chi)\otimes\u^-(\chi)$ such that
\begin{equation}\label{eqn:def eta}
\eta(K_i,L_j)=q_{ij}, \qquad \eta(E_i,F_j)=-\delta_{ij}, \qquad \eta(E_i,L_j)=\eta(K_i,F_j)=0.
\end{equation}
for all $1\leq i,j\leq\theta$. It satisfies the following condition: for all $E\in\u^+(\chi)$, $F\in\u^-(\chi)$, $K\in\u^{+0}(\chi)$, $L\in\u^{-0}(\chi)$
\begin{equation}\label{eqn:decomp eta}
\eta(EK,FL)=\eta(E,F)\eta(K,L).
\end{equation}
Moreover, if $\beta\neq\gamma\in\N_0^\theta$, then $\eta|_{\u^+(\chi)_\beta\otimes\u^-(\chi)_{-\gamma}} \cong 0$.\qed
\end{prop}

\bigbreak

Assume that all the integers $a_{ij}^\upsilon$
are defined, $\upsilon\in\cG(\chi)$, so the automorphisms
$s_p^\upsilon$ are defined. For simplicity, we denote $\Eb_i$,
$\Fb_i$, $\Kb_i$, $\Lb_i$ the generators corresponding to
$\cU(\schi)$, $a_{ij}=a_{ij}^\chi$, $q_{ij}=q_{ij}^\chi$,
$\qb_{ij}=q_{ij}^\schi$. We define also the scalars
\begin{equation}\label{eqn:escalares lambda}
\lambda_i(\chi):=
(-a_{pi})_{q_{pp}}\prod_{s=0}^{-a_{pi}-1}(q_{pp}^sq_{pi}q_{ip}-1) ,
\qquad i \neq p.
\end{equation}

Fix $p\in\{1,\ldots,\theta\}$. If $i\neq p$ we consider the elements
\cite{H-isom},
$$ E_{i,0(p)}^+, E_{i,0(p)}^-:=E_i, \qquad F_{i,0(p)}^+, F_{i,0(p)}^-:=F_i, $$
and recursively,
\begin{align*}
E_{i,m+1(p)}^+ &:= E_pE_{i,m(p)}^+ - (K_p \cdot E_{i,m(p)}^+)E_p =
(\ad_c E_p)^{m+1}E_i,
\\ E_{i,m+1(p)}^- &:= E_pE_{i,m(p)}^- - (L_p \cdot E_{i,m(p)}^-)E_p,
\\ F_{i,m+1(p)}^+ &:= F_pF_{i,m(p)}^+ - (L_p \cdot F_{i,m(p)}^+)F_p,
\\ F_{i,m+1(p)}^- &:= F_pF_{i,m(p)}^- - (K_p \cdot F_{i,m(p)}^-)F_p.
\end{align*}
If $p$ is explicit, we simply denote $E_{i,m(p)}^\pm$ by
$E_{i,m}^\pm$. By \cite[Corollary 5.4]{H-isom},
\begin{equation}\label{eqn:corchete E+ con Fp}
E_{i,m}^+ F_i - F_i E_{i,m}^+=
(m)_{q_{pp}}(q_{pp}^{m-1}q_{pi}q_{ip}-1)L_p E_{i,m-1}^+.
\end{equation}

\begin{theorem}
There exist algebra morphisms
\begin{equation}\label{eqn:iso lusztig para cUp}
T_p, T_p^-: \u(\chi) \to \u(\schi)
\end{equation}
univocally determined by the following conditions:
\begin{align*}
T_p(K_p)&=T_p^-(K_p)=\Kb_p^{-1}, & T_p(K_i)&=T_p^-(K_i)=\Kb_p^{m_{pi}}\Kb_i,
\\ T_p(L_p)&=T_p^-(L_p)=\Lb_p^{-1}, & T_p(L_i)&=T_p^-(L_i)=\Lb_p^{m_{pi}}\Lb_i,
\\ T_p(E_p)&=\Fb_p\Lb_p^{-1}, & T_p(E_i)&=\Eb^+_{i,m_{pi}},
\\ T_p(F_p)&=\Kb_p^{-1}\Eb_p, & T_p(F_i)&=\lambda_p(\schi)^{-1}\Fb^+_{i,m_{pi}},
\\ T_p^-(E_p)&=\Kb_p^{-1}\Fb_p, & T_p^-(E_i)&=\lambda_p(\schi^{-1})^{-1}\Eb^-_{i,m_{pi}},
\\ T_p^-(F_p)&=\Eb_p\Lb_p^{-1}, & T_p^-(F_i)&=\Fb^-_{i,m_{pi}}.
\end{align*}
for every $i\neq p$. Moreover, $T_pT_p^-=T_p^-T_p=\id$, and there
exists $\mu \in (\ku^\times)^\theta$ such that
\begin{equation}\label{eqn:Tp con phi4}
T_p \circ \phi_4 = \phi_4 \circ T_p^- \circ \varphi_{\mu}.
\end{equation}
\qed
\end{theorem}

By \cite[Proposition 4.2]{HY-shapov}, we have for all $\alpha\in\Z^\theta$
\begin{align}\label{eq:degree applying Tp}
T_p\big(\u(\chi)_\alpha \big)= \u(\schi)_{s_p^\chi(\alpha)}.
\end{align}

\section{$R$-matrix from a version of a universal $R$-matrix} \label{section:Rfnrep}

Most of the ideas we shall give in this section are modifications of \cite[Section~4]{Tanisaki}.
Let $\chi : \Zt \times \Zt \to \ku^{\times}$ be any bicharacter. We will compute an
R-matrix for some modules of $\truchi$ from canonical elements of $\truchi$.
If $M=|\Delta_+^\chi|<\infty$, the canonical elements can be obtained by
Proposition \ref{prop:duality PBW bases}.

\subsection{Equations for canonical elements}\label{subsection:Efce}

We recall \cite[(3.18), (3.19)]{HY-shapov}:
\begin{align}
YX&=\treta(X_{1},\cS(Y_{1}))\treta(X_3,Y_3)X_2 Y_2, \label{eqn:DbleComYX} \\
XY&=\treta(X_1,Y_1)\treta(X_3,\cS(Y_3))Y_2X_2, & X&\in\trugeqchi, \, Y\in\truleqchi.  \label{eqn:DbleComXY}
\end{align}

Define the $\ku$-linear homomorphism $\tau:\truchi\otimes\truchi\to\truchi\otimes\truchi$ by
$$ \tau(X\otimes Y):=Y\otimes X.$$
Given $\textbf{X}\in\trugeqchi$, $\textbf{Y}\in\truleqchi$, we define the $\ku$-linear homomorphisms
\begin{align*}
\hattretaleq_{\textbf{X}}&:\truleqchi\to\ku, & \hattretaleq_{\textbf{X}}(Y)&:=\treta(\textbf{X},Y), & Y\in&\truleqchi, \\
\hattretageq_{\textbf{Y}}&:\trugeqchi\to\ku, & \hattretageq_{\textbf{Y}}(X)&:=\treta(X,\textbf{Y}), & X\in&\trugeqchi.
\end{align*}

\begin{lemma} \label{lemma:preCom}
Let $1\leq i\leq\trtheta$ and $\trbeta\in\Nzerotheta$. Set
$$ \Nzerotheta(\trbeta;i):=\{\,\trgamma\in\Nzerotheta\ -\{0,\tral_i\}, |\,\trbeta-\trgamma\in\Nzerotheta -\{0,\tral_i\}\}. $$

\noindent \vi Let $\trbeta\notin\{0,\tral_i,2\tral_i\}$, $Y\in\trumchi_{-\trbeta}$.
Set $Y^\prime$, $Y^{\prime\prime}\in\trumchi_{-\trbeta+\tral_i}$ such that
$[\trE_i,Y]=\trK_iY^\prime-Y^{\prime\prime}\trL_i$. Then
\begin{align} \label{eqn:FiYmYFi}
\trDelta(Y)- &(Y\otimes\trL^\trbeta+1\otimes Y
+\trF_i\otimes Y^{\prime\prime}\trL^{\tral_i}
+Y^\prime\otimes\trF_i\trL^{\trbeta-\tral_i})
\\ &\in \oplus_{\trgamma\in\Nzerotheta(\trbeta;i)}
\trumchi_{-\trgamma}\otimes\trumchi_{-\trbeta+\trgamma}\trL^\trgamma.\nonumber
\end{align}
In particular,
\begin{align} \label{eqn:FiYmYFid}
(\hattretaleq_{\trE_i}\otimes\rmid)(\trDelta(Y))&=-Y^{\prime\prime}\trL_i, & (\rmid\otimes\hattretaleq_{\trE_i})(\trDelta(Y))&=-Y^\prime.
\end{align}

\noindent \vii Let $\trbeta\notin\{0,\tral_i,2\tral_i\}$, $X\in\trupchi_\trbeta$.
Set $X^\prime, X^{\prime\prime}\in\trupchi_{\trbeta-\tral_i}$ such that
$[X,\trF_i]=X^{\prime\prime}\trK_i-\trL_iX^\prime$. Then
\begin{align} \label{eqn:EiXmXEi}
\trDelta(X)- &(X\otimes 1+\trK^\trbeta \otimes X
+X^{\prime\prime}\trK^{\tral_i} \otimes\trE_i
+\trE_i\trK^{\trbeta-\tral_i}\otimes X^\prime)
\\ &\in\oplus_{\trgamma\in\Nzerotheta(\trbeta;i)}
\trupchi_\trgamma\trK^{\trbeta-\trgamma}\otimes\trupchi_{\trbeta-\trgamma}. \nonumber
\end{align}
In particular,
\begin{align} \label{eqn:EiXmXEid}
(\rmid\otimes\hattretageq_{\trF_i})(\trDelta(X))&=-X^{\prime\prime}\trK_i, & (\hattretageq_{\trF_i}\otimes\rmid)(\trDelta(X))&=-X^\prime.
\end{align}
\end{lemma}
\pf
We prove \emph{\vi}; \emph{\vii} can be proved analogously. Note that
$$\trDeltatwo(\trE_i)=\trE_i\otimes 1\otimes 1+ \trK_i\otimes\trE_i\otimes 1
+\trK_i\otimes\trK_i\otimes\trE_i.$$
Define ${\bar Y}^\prime$, ${\bar Y}^{\prime\prime}$
as the elements of $\trumchi_{-\trbeta+\tral_i}$
satisfying the same property as \eqref{eqn:FiYmYFi}
with ${\bar Y}^\prime$, ${\bar Y}^{\prime\prime}$
in place of $Y^\prime$, $Y^{\prime\prime}$.
By \eqref{eqn:DbleComYX}, we have
\begin{align*}
Y\trE_i = &\treta(\trE_i,\trS(\trF_i))\treta(1,\trL^\trbeta){\bar Y}^{\prime\prime}\trL^{\tral_i}
+\treta(\trK_i,\trS(1))\treta(1,\trL^\trbeta)\trE_iY
\\ &\quad +\treta(\trK_i,\trS(1))\treta(\trE_i,\trF_i\trL^{\trbeta-\tral_i})\trK_i{\bar Y}^\prime
\\ =&\treta(\trE_i,-\trF_i\trL^{-\tral_i})\treta(1,\trL^\trbeta){\bar Y}^{\prime\prime}\trL^{\tral_i}
+\treta(\trK_i,1)\treta(1,\trL^\trbeta)\trE_iY
\\ &\quad +\treta(\trK_i,1)\treta(\trE_i,\trF_i\trL^{\trbeta-\tral_i})\trK_i{\bar Y}^\prime
\\ = &{\bar Y}^{\prime\prime}\trL_i +\trE_iY -\trK_i{\bar Y}^\prime,
\end{align*}
so the proof is complete.
\epf

\smallskip

Fix $\trbeta\in\Nzerotheta$ and $m_{\beta}:=\dim \trupchi_\trbeta=\dim \trumchi_{-\trbeta}$.
Fix also $\{\trE^{(\trbeta)}_x\}$, $\{\trF^{(\trbeta)}_y\}$ bases of the spaces $\trupchi_\trbeta$, $\trumchi_{-\trbeta}$, which are dual for $\eta$. The matrix $[\treta(\trE^{(\trbeta)}_x,\trF^{(\trbeta)}_y)]_{1\leq x,y\leq m_{\beta}}$ is invertible,
we call $[b^{(\trbeta)}_{xy}]_{1\leq x,y\leq m_{\beta}}$ to its inverse.

\begin{lemma}\label{lemma:dleta} For all $X\in\trupchi_\trbeta$, $ Y\in\trumchi_{-\trbeta}$ it holds:
\begin{align}
X=& \sum_{x,y}b^{(\trbeta)}_{yx}\treta(X,\trF^{(\trbeta)}_y)\trE^{(\trbeta)}_x
\label{eqn:dletaX} \\
X=& \sum_{x,y}b^{(\trbeta)}_{yx}\treta(\trE^{(\trbeta)}_x,Y)\trF^{(\trbeta)}_y. \label{eqn:dletaY}
\end{align}
\end{lemma}
\pf
We prove \eqref{eqn:dletaX}; the proof of \eqref{eqn:dletaY} is similar. We have
\begin{align*}
\treta\left(\sum_{x,y}b^{(\trbeta)}_{yx}\treta(X,\trF^{(\trbeta)}_y) \trE^{(\trbeta)}_x,\trF^{(\trbeta)}_z\right) &=\sum_{x,y}b^{(\trbeta)}_{yx}\treta(X,\trF^{(\trbeta)}_y)
\treta(\trE^{(\trbeta)}_x,\trF^{(\trbeta)}_z) \\
&=\sum_y\delta_{yz}\treta(X,\trF^{(\trbeta)}_y)=\treta(X,\trF^{(\trbeta)}_z),
\end{align*}
for all $1\leq z\leq m$. \eqref{eqn:dletaX} follows since $\treta_{|\trupchi_\trbeta\times\trumchi_{-\trbeta}}$
is non-degenerate.
\epf

\smallskip

Let $\mcC_\trbeta$ be the canonical element of $\trupchi_\trbeta\otimes\trumchi_{-\trbeta}$, i.e.
\begin{equation*}
\mcC_\trbeta=\sum_{x,y=1}^{m_{\beta}}
b^{(\trbeta)}_{yx}\trE^{(\trbeta)}_x\otimes\trF^{(\trbeta)}_y.
\end{equation*}

\begin{lemma}\label{lemma:Com} Let $1\leq i\leq\trtheta$. The following identities hold:
\begin{align}
[1\otimes\trE_i,\mcC_{\trbeta+\tral_i}] &= \mcC_\trbeta(\trE_i\otimes\trL_i)-(\trE_i\otimes\trK_i)\mcC_\trbeta, \label{eqn:ECCE} \\
[\mcC_{\trbeta+\tral_i},F_i\otimes 1 ] &= (\trL_i\otimes\trF_i)\mcC_\trbeta-\mcC_\trbeta(\trK_i\otimes\trF_i). \label{eqn:FCCF}
\end{align}
\end{lemma}
\pf
We prove \eqref{eqn:ECCE}.
Let $Y\in\trumchi_{-\trbeta-\tral_i}$.
Let $Y^\prime$, $Y^{\prime\prime}\in\trumchi_{-\trbeta}$
be such that $[\trE_i,Y]=Y^\prime\trK_i-\trL_iY^{\prime\prime}$. Using \eqref{eqn:dletaX}, we have
\begin{align}
(\hattretageq_Y\otimes\rmid)([1\otimes\trE_i,\mcC_{\trbeta+\tral_i}]) & =\sum_{x,y}b^{(\trbeta+\tral_i)}_{yx} \treta(\trE^{(\trbeta+\tral_i)}_x,Y) [\trE_i,\trF^{(-\trbeta-\tral_i)}_y] \nonumber \\
\label{eqn:ECCEa}  = & \left[\trE_i,\sum_{x,y}b^{(\trbeta+\tral_i)}_{yx} \treta(\trE^{(\trbeta+\tral_i)}_x,Y)\trF^{(-\trbeta-\tral_i)}_y\right] =[\trE_i,Y].
\end{align}
Now using \eqref{eqn:FiYmYFid}, \eqref{eqn:dletaX} and \eqref{eqn:ECCEa}, we compute:
\begin{align*}
(\hattretageq_Y\otimes & \rmid)(\mcC_\trbeta(\trE_i\otimes\trL_i)-(\trE_i\otimes\trK_i)\mcC_\trbeta) \\
&=\sum_{x,y}b^{(\trbeta)}_{yx} (\treta(\trE^{(\trbeta)}_x\trE_i,Y)\trF^{(\trbeta)}_y\trL_i- \treta(\trE_i\trE^{(\trbeta)}_x,Y)\trK_i\trF^{(\trbeta)}_y) \\
&=\sum_{x,y}b^{(\trbeta)}_{yx} (\treta(\trE_i\otimes\trE^{(\trbeta)}_x,\trDelta(Y))\trF^{(\trbeta)}_y\trL_i-\treta(\trE^{(\trbeta)}_x\otimes\trE_i,\trDelta(Y))\trK_i\trF^{(\trbeta)}_y) \\
&=\sum_{x,y}b^{(\trbeta)}_{yx} (-\treta(\trE^{(\trbeta)}_x,Y^{\prime\prime})\trF^{(\trbeta)}_y\trL_i+ \treta(\trE^{(\trbeta)}_x,Y^\prime)\trK_i\trF^{(\trbeta)}_y)\\
&=-Y^{\prime\prime}\trL_i+\trK_iY^\prime =[\trE_i,Y] =(\hattretageq_Y\otimes\rmid)([1\otimes\trE_i,\mcC_{\trbeta+\tral_i}]).
\end{align*}
Since $\treta_{|\trupchi_\trbeta\times\trumchi_{-\trbeta}}$ is non-degenerate, we have \eqref{eqn:ECCE}. Similarly we obtain \eqref{eqn:FCCF}.
\epf

\smallskip

\begin{lemma} \label{lemma:preCinv}
Let $\mcC^\prime_\trbeta:=(\trK^\trbeta\otimes 1)(\cS\otimes\rmid)(\mcC_\trbeta)$. For every $\tral\in\Nzerotheta$,
\begin{equation}\label{eqn:CCp}
\sum_{{{\trbeta,\,\trgamma\in\Nzerotheta}
\atop{\trbeta+\trgamma=\tral}}}\mcC_\trbeta\mcC^\prime_\trgamma
=\delta_{\alpha,0}=\sum_{{{\trbeta,\,\trgamma\in\Nzerotheta}
\atop{\trbeta+\trgamma=\tral}}}\mcC^\prime_\trbeta\mcC_\trgamma.
\end{equation}
\end{lemma}
\pf
If $\tral=0$, \eqref{eqn:CCp} is clear. Assume $\tral\ne0$.
We show the first equation of \eqref{eqn:CCp}.
Since
$\treta_{|\trupchi_\tral\times\trumchi_{-\tral}}$ is non-degenerate,
it suffices to show that
\begin{align}\label{eqn:ztreta}
\sum_{{{\trbeta,\,\trgamma\in\Nzerotheta} \atop{\trbeta+\trgamma=\tral}}}(\hattretageq_Y\otimes\rmid_\truchi)(\mcC_\trbeta\mcC^\prime_\trgamma)&=0,   &
\mbox{for all }&Y\in\trumchi_{-\tral}.
\end{align}
Write $\trDelta(Y)=\sum_{{{\trbeta,\,\trgamma\in\Nzerotheta}
\atop{\trbeta+\trgamma=\tral}}}Y^{(\trbeta,\trgamma)}(1\otimes\trL^\trbeta )$,
where $Y^{(\trbeta,\trgamma)}\in\trumchi_{-\trbeta}\otimes\trumchi_{-\trgamma}$.
Further write $Y^{(\trbeta,\trgamma)}
=\sum_mY^{(\trbeta,\trgamma)}_{-\trbeta,m}\otimes Y^{(\trbeta,\trgamma)\prime}_{-\trgamma,m}$,
where $Y^{(\trbeta,\trgamma)}_{-\trbeta,m}\in\trumchi_{-\trbeta}$
and $Y^{(\trbeta,\trgamma)\prime}_{-\trgamma,m}\in\trumchi_{-\trgamma}$. The left hand side of \eqref{eqn:ztreta} is
\begin{align*}
\sum_{{{\trbeta,\,\trgamma\in\Nzerotheta} \atop{\trbeta+\trgamma=\tral}}} \sum_{{{x,x^\prime,} \atop{y,y^\prime}}}
b^{(\trbeta )}_{yx}b^{(\trgamma )}_{y^\prime x^\prime} & \treta (\trE^{(\trbeta )}_x\trK^\trgamma\cS(\trE^{(\trgamma )}_{x^\prime}),Y) \trF^{(\trbeta )}_y\trF^{(\trgamma )}_{y^\prime}
\\
=\sum_{{{\trbeta,\,\trgamma\in\Nzerotheta} \atop{\trbeta+\trgamma=\tral}}} \sum_{{{m,x,y,} \atop{x^\prime,y^\prime}}} &
b^{(\trbeta )}_{yx}b^{(\trgamma )}_{y^\prime x^\prime}
\treta (\trE^{(\trbeta )}_x,Y^{(\trbeta,\trgamma)\prime}_{\trgamma,m}\trL^\trbeta )
\treta (\trK^\trgamma\cS(\trE^{(\trgamma )}_{x^\prime}),Y^{(\trbeta,\trgamma)}_{\trbeta,m})
\trF^{(\trbeta )}_y\trF^{(\trgamma )}_{y^\prime}
\\
=\sum_{{{\trbeta,\,\trgamma\in\Nzerotheta} \atop{\trbeta+\trgamma=\tral}}} \sum_{{{m,x,y,} \atop{x^\prime,y^\prime}}} &
b^{(\trbeta )}_{yx}b^{(\trgamma )}_{y^\prime x^\prime}
\treta (\trE^{(\trbeta )}_x,Y^{(\trbeta,\trgamma)\prime}_{\trgamma,m}\trL^\trbeta )
\treta (\cS(\trE^{(\trgamma )}_{x^\prime}\trK^{\trgamma} ),Y^{(\trbeta,\trgamma)}_{\trbeta,m})
\trF^{(\trbeta )}_y\trF^{(\trgamma )}_{y^\prime}
\\
=\sum_{{{\trbeta,\,\trgamma\in\Nzerotheta} \atop{\trbeta+\trgamma=\tral}}} \sum_{{{m,x,y,} \atop{x^\prime,y^\prime}}} &
b^{(\trbeta )}_{yx}b^{(\trgamma )}_{y^\prime x^\prime}
\treta (\trE^{(\trbeta )}_x,Y^{(\trbeta,\trgamma)\prime}_{\trgamma,m}\trL^\trbeta )
\treta (\trE^{(\trgamma )}_{x^\prime}\trK^{\trgamma} ,\cS^{-1}(Y^{(\trbeta,\trgamma)}_{\trbeta,m}))
\trF^{(\trbeta )}_y\trF^{(\trgamma )}_{y^\prime}
\\
=\sum_{{{\trbeta,\,\trgamma\in\Nzerotheta} \atop{\trbeta+\trgamma=\tral}}} \sum_{{{m,x,y,} \atop{x^\prime,y^\prime}}} &
b^{(\trbeta )}_{yx}b^{(\trgamma )}_{y^\prime x^\prime}
\treta (\trE^{(\trbeta )}_x,Y^{(\trbeta,\trgamma)\prime}_{\trgamma,m}\trL^\trbeta )
\treta (\trE^{(\trgamma )}_{x^\prime}\trK^{\trgamma} ,\cS^{-1}(Y^{(\trbeta,\trgamma)}_{\trbeta,m})) \trL^\trgamma
\\
=\sum_{{{\trbeta,\,\trgamma\in\Nzerotheta} \atop{\trbeta+\trgamma=\tral}}} \sum_{{{x,x^\prime,} \atop{y,y^\prime}}}
b^{(\trgamma )}_{y^\prime x^\prime}&b^{(\trbeta )}_{yx} \treta(\trK^{\trgamma},\trL^{\trgamma}) ( \sum_m \treta (\trE^{(\trbeta )}_x,Y^{(\trbeta,\trgamma)\prime}_{\trgamma,m} ) \treta (\trE^{(\trgamma )}_{x^\prime},\cS^{-1}(Y^{(\trbeta,\trgamma)}_{\trbeta,m}) )\trL^\trgamma
\\
=\sum_{{{\trbeta,\,\trgamma\in\Nzerotheta} \atop{\trbeta+\trgamma=\tral}}} \sum_m & \trchi(\trgamma,\trgamma) Y^{(\trbeta,\trgamma)\prime}_{\trgamma,m}\cS^{-1}(Y^{(\trbeta,\trgamma)}_{\trbeta,m})\trL^\trgamma
\\
=\sum_{{{\trbeta,\,\trgamma\in\Nzerotheta} \atop{\trbeta+\trgamma=\tral}}} \sum_m & Y^{(\trbeta,\trgamma)\prime}_{\trgamma,m}\trL^\trgamma \cS^{-1}(Y^{(\trbeta,\trgamma)}_{\trbeta,m}) =\varepsilon (Y)=0,
\end{align*}
where we use \eqref{eqn:dletaY} and the graduation of $\u(\chi)$. The second equation of \eqref{eqn:CCp} is obtained in a similar way.
\epf

\medbreak

\begin{lemma} \label{eqn:idDel}
The following identities hold:
\begin{align}
& (\rmid\otimes\trDelta)(\mcC_\tral)=\sum_{\trbeta+\trgamma=\tral}\mcC_\trbeta^{(1,3)}\mcC_\trgamma^{(1,2)}
(1\otimes 1\otimes \trL^\trgamma), \label{eqn:idotDel}
 \\
& (\trDelta\otimes\rmid)(\mcC_\tral)=\sum_{\trbeta+\trgamma=\tral}\mcC_\trbeta^{(1,3)}\mcC_\trgamma^{(2,3)}
(\trK^\trgamma\otimes 1\otimes 1).
\label{eqn:Delotid}
\end{align}
\end{lemma}
\pf
We show \eqref{eqn:idotDel}. Given $X_1\in\trupchi_\trgamma$ and $X_2\in\trupchi_\trbeta$, we compute
\begin{align*}
(\rmid\otimes\hattretaleq_{X_1}&\otimes\hattretaleq_{X_2})(\rmid\otimes\trDelta)(\mcC_\tral)= \sum_{x,y}b^{(\tral)}_{yx}\treta(X_2X_1,\trF^{(-\tral)}_y)\trE^{(\tral)}_x  \\
&=X_2X_1 = \sum_{x^{\prime\prime},y^{\prime\prime},x^\prime,y^\prime} b^{(\trbeta)}_{y^{\prime\prime} x^{\prime\prime}}b^{(\trgamma)}_{y^\prime x^\prime}
\treta(X_2,\trF^{(\trbeta)}_{y^{\prime\prime}})\treta(X_1,\trF^{(-\trgamma)}_{y^\prime})
\trE^{(\trbeta)}_{x^{\prime\prime}}\trE^{(\trgamma)}_{x^\prime} \\
&= (\rmid\otimes\hattretaleq_{X_1}\otimes\hattretaleq_{X_2})(\mcC_\trbeta^{(1,3)}\mcC_\trgamma^{(1,2)}),
\end{align*}
where we use \eqref{eqn:dletaX} twice. Since
$$ (\rmid\otimes\trDelta)(\mcC_\tral)\in\sum_{\trbeta+\trgamma=\tral} \trupchi_\tral\otimes\trupchi_\trgamma\otimes\trupchi_\trbeta\trL^\trgamma,$$
we prove that \eqref{eqn:idotDel} holds. Similarly we obtain \eqref{eqn:Delotid}.
\epf

\subsection{$R$-matrix for finite dimensional $\truchi$-modules}\label{subsection:Rfinmod}
Fix $\trV_1$, $\trV_2$, $\trV_3$ three finite dimensional $\truchi$-modules, with associated $\ku$-algebra homomorphisms
$\trrho_x:\truchi\to\rmEnd_\ku(\trV_x)$, $x\in\{1,2,3\}$, such that there exist an element $\trv_x\in\trV_x$ and a $\ku$-algebra homomorphism $\trLambda_x:\truzerochi\to\ku$ for each $x\in\{1,2,3\}$ satisfying
\begin{align*}
X\cdot\trv_x &=\trLambda_x(X)\trv_x\mbox{ for all }X\in\truzerochi, & \trV_x&=\trumchi\cdot\trv_x, \\
\trE_i\cdot\trv_x&=0\mbox{ for all }1\leq i\leq \theta.
\end{align*}
If $\mcF=\sum_z\mcF^\prime_z\otimes\mcF^{\prime\prime}_z \in\rmEnd_\ku(\trV_x\otimes\trV_y)
\cong\rmEnd_\ku(\trV_x)\otimes\rmEnd_\ku(\trV_y)$, $1\leq x<y\leq 3$, we set $\mcF^{(x,y)} \in \rmEnd_\ku(\trV_1\otimes\trV_2\otimes\trV_3)
\cong\rmEnd_\ku(\trV_1)\otimes\rmEnd_\ku(\trV_2)\otimes\rmEnd_\ku(\trV_3)$ as
\begin{align*}
\mcF^{(x,y)}&=\sum_z\mcF^\prime_z\otimes\mcF^{\prime\prime}_z\otimes\rmid_{\trV_3} & \mbox{if }&x=1, y=2, \\
\mcF^{(x,y)}&=\sum_z\mcF^\prime_z\otimes\rmid_{\trV_2}\otimes\mcF^{\prime\prime}_z & \mbox{if }&x=1, y=3, \\
\mcF^{(x,y)}&=\rmid_{\trV_1}\otimes\sum_z\mcF^\prime_z\otimes\mcF^{\prime\prime}_z & \mbox{if }&x=2, y=3.
\end{align*}
Now define $\trf_{xy}\in\rmGL_\ku(\trV_x\otimes\trV_y)$ by
\begin{equation*}
\trf_{xy}(X\trv_x\otimes Y\trv_y):=\trchi(\trbeta,\tral)\trLambda_x(\trK^{-\trbeta})\trLambda_y(\trL^\tral)X\trv_x\otimes Y\trv_y
\end{equation*}
for $\tral$, $\trbeta\in\Nzerotheta$ and $X\in\trumchi_{-\tral}$, $Y\in\trumchi_{-\trbeta}$. Set also
\begin{align*}
\cC_{xy}&:=\sum_{\trbeta\in\Nzerotheta}(\trrho_x\otimes\trrho_y)(\mcC_\trbeta), & \trR_{xy}&:=\cC_{xy}\trf_{xy}^{-1}.
\end{align*}

\begin{lemma} \label{lemma:premAPP} For each $1\leq i\leq\trtheta$ and $\cheX\in\trV_x\otimes\trV_y$,
\begin{align}
\trf_{xy}((\trE_i\otimes 1)\cheX)&=(\trE_i\otimes\trL_i^{-1})\trf_{xy}(\cheX), \label{eqn:fxyEoto} \\
\trf_{xy}((1\otimes\trE_i)\cheX)&=(\trK_i\otimes\trE_i)\trf_{xy}(\cheX), \label{eqn:fxyotoE} \\
\trf_{xy}((\trF_i\otimes 1)\cheX)&=(\trF_i\otimes\trL_i)\trf_{xy}(\cheX), \label{eqn:fxyFoto} \\
\trf_{xy}((1\otimes\trF_i)\cheX)&=(\trK_i^{-1}\otimes\trF_i)\trf_{xy}(\cheX).
\label{eqn:fxyotoF}
\end{align}

\end{lemma}
\pf
We show \eqref{eqn:fxyEoto}. For each $X\in\trumchi_{-\trbeta}$, $Y\in\trumchi_{-\trgamma}$,
\begin{align*}
\trf_{xy}((\trE_i\otimes 1)X\trv_x\otimes Y\trv_y) & =\trf_{xy}(\trE_iX\trv_x\otimes Y\trv_y) \\
&=\trchi(\trgamma,\trbeta-\tral_i)\trLambda_x(\trK^{-\trgamma})\trLambda_y(\trL^{\trbeta-\tral_i}) \trE_iX\trv_x\otimes Y\trv_y \\
&=(\trE_i\otimes \trL_i^{-1})\trf_{xy}(X\trv_x\otimes Y\trv_y).
\end{align*}
Thus we have \eqref{eqn:fxyEoto}.
Similarly we obtain \eqref{eqn:fxyotoE},
\eqref{eqn:fxyFoto}
and \eqref{eqn:fxyotoF}.
\epf

Now we are ready to obtain the $R$-matrix for the modules $V_x$, $1\leq x\leq3$.

\begin{theorem} \label{theorem:mainAPP}
\vi $\cC_{xy}\in\rmGL_\ku(\trV_x\otimes\trV_y)$ and
\begin{equation}\label{eqn:invCxy}
 \cC_{xy}^{-1}=\sum_{\trbeta\in\Nzerotheta}(\trrho_x\otimes\trrho_y)(\trK^\trbeta\otimes 1)(\cS\otimes\rmid)(\mcC_\beta).
\end{equation}
\vii For every $X\in\truchi$
\begin{equation}\label{eqn:RCom}
\trR_{xy}(\trrho_x\otimes\trrho_y)(\trDelta(X))\trR_{xy}^{-1}
=(\trrho_x\otimes\trrho_y)((\trtau\circ\trDelta)(X)).
\end{equation}
\viii The following identities hold:
\begin{align}
\sum_{\trbeta\in\Nzerotheta}(\trrho_1\otimes\trrho_2\otimes\trrho_3)((\trDelta\otimes\rmid_\truchi)(\mcC_\beta))
&=\cC_{13}^{(1,3)}(\trf_{13}^{(1,3)})^{-1}\cC_{23}^{(2,3)}\trf_{13}^{(1,3)},
\label{eqn:preRacbc}\\
\sum_{\trbeta\in\Nzerotheta}(\trrho_1\otimes\trrho_2\otimes\trrho_3)((\rmid_\truchi\otimes\trDelta)(\mcC_\beta))
&=\cC_{13}^{(1,3)}(\trf_{13}^{(1,3)})^{-1}\cC_{12}^{(1,2)}\trf_{13}^{(1,3)}.
\label{eqn:preRacab}
\end{align}
\viv The elements $R_{xy}$ satisfy:
\begin{equation}
\trR_{12}^{(1,2)}\trR_{13}^{(1,3)}\trR_{23}^{(2,3)}
=\trR_{23}^{(2,3)}\trR_{13}^{(1,3)}\trR_{12}^{(1,2)}.
\end{equation}
\end{theorem}
\pf
\emph{\vi} This immediately follows from \eqref{eqn:CCp}.
\smallskip

\noindent \emph{\vii} As we have algebra maps on both sides of the identity, it is enough to prove it for the generators of $\truchi$, and it follows by using Lemmata \ref{lemma:Com}, \ref{lemma:premAPP}.
For example, for each $\cheX\in V_x\otimes V_y$, by \eqref{eqn:fxyEoto}, \eqref{eqn:fxyotoE}, \eqref{eqn:ECCE} we have
\begin{align*}
(\trR_{xy}\trDelta & (\trE_i)-(\trtau\circ\trDelta)(\trE_i)\trR_{xy})\cheX =(\cC_{xy}\trf_{xy}^{-1}\trDelta(\trE_i)-(\trtau\circ\trDelta)(\trE_i)\cC_{xy}\trf_{xy}^{-1})\cheX \\
&=\sum_{\trbeta\in\Nzerotheta} (\mcC_\trbeta\trf_{xy}^{-1}(\trE_i\otimes 1+\trK_i\otimes\trE_i)-(1\otimes \trE_i+\trE_i\otimes\trK_i)\mcC_\trbeta\trf_{xy}^{-1})\cheX \\
&=\sum_{\trbeta\in\Nzerotheta} (\mcC_\trbeta(\trE_i\otimes \trL_i+1\otimes\trE_i)-(1\otimes \trE_i+\trE_i\otimes\trK_i)\mcC_\trbeta)\trf_{xy}^{-1}\cheX \\
&=\sum_{\trbeta\in\Nzerotheta}([1\otimes\trE_i,\mcC_{\trbeta+\tral_i}]-[1\otimes\trE_i,\mcC_\trbeta])\trf_{xy}^{-1}\cheX  \\
&= -\sum_{\trbeta\in\Nzerotheta,\trbeta-\tral_i\notin\Nzerotheta} [1\otimes\trE_i,\mcC_\trbeta]\trf_{xy}^{-1}\cheX=0.
\end{align*}

\noindent \emph{\viii} It can be proved by using Lemmata \ref{eqn:idDel}, \ref{lemma:premAPP}. In fact, we compute for each $\cheX\in V_x\otimes V_y$:
\begin{align*}
\cC_{13}^{(1,3)}(\trf_{13}^{(1,3)})^{-1}\cC_{23}^{(2,3)}\trf_{13}^{(1,3)}\cheX
&= \sum_{\tral,\,\trgamma\in\Nzerotheta}\mcC^{(1,3)}_\tral(\trf_{13}^{(1,3)})^{-1}(\mcC^{(2,3)}_\trgamma \trf_{13}^{(1,3)}(\cheX)) \\
&=\sum_{\tral,\,\trgamma\in\Nzerotheta}\mcC^{(1,3)}_\tral\mcC^{(2,3)}_\trgamma (\trK^\trgamma\otimes 1\otimes 1)\cheX \\
&= \sum_{\trbeta\in\Nzerotheta}(\trrho_1\otimes\trrho_2\otimes\trrho_3)((\trDelta\otimes\rmid_\truchi)(\mcC_\beta))\cheX.
\end{align*}

\noindent \emph{\viv}  In this case the proof follows by Lemma \ref{lemma:Com} and the previous claims:
\begin{align*}
\trR_{12}^{(1,2)}&\trR_{13}^{(1,3)}\trR_{23}^{(2,3)} =\trR_{12}^{(1,2)}\cC_{13}^{(1,3)}(\trf_{13}^{(1,3)})^{-1}\cC_{23}^{(2,3)}(\trf_{23}^{(2,3)})^{-1} \\
&=\sum_{\trbeta\in\Nzerotheta}\trR_{12}^{(1,2)} (\trrho_1\otimes\trrho_2\otimes\trrho_3)((\trDelta\otimes\rmid_\truchi)(\mcC_\beta))
(\trf_{13}^{(1,3)})^{-1}(\trf_{23}^{(2,3)})^{-1}\\
&=\sum_{\trbeta\in\Nzerotheta} (\trrho_1\otimes\trrho_2\otimes\trrho_3)((\trtau\circ\trDelta)\otimes\rmid_\truchi)(\mcC_\beta)
\trR_{12}^{(1,2)} (\trf_{13}^{(1,3)})^{-1}(\trf_{23}^{(2,3)})^{-1}\\
&=\cC_{23}^{(2,3)}(\trf_{23}^{(2,3)})^{-1}\cC_{13}^{(1,3)}\trf_{23}^{(2,3)}\trR_{12}^{(1,2)}
(\trf_{13}^{(1,3)})^{-1}(\trf_{23}^{(2,3)})^{-1}\\
&=\trR_{23}^{(2,3)}\trR_{13}^{(1,3)}\trf_{13}^{(1,3)}\trf_{23}^{(2,3)}\trR_{12}^{(1,2)}(\trf_{13}^{(1,3)})^{-1}(\trf_{23}^{(2,3)})^{-1} \\
&=\trR_{23}^{(2,3)}\trR_{13}^{(1,3)}\trR_{12}^{(1,2)}.
\end{align*}
\epf

\section{R-matrices of quantum doubles of Nichols algebras with finite root systems}\label{section:universal Rmatrix}

For this section we fix $\chi$ such that $M=|\Delta_+^\chi|<\infty$. First we recall a series of
results from \cite[Section 4]{HY-shapov}, which will be useful to compute explicitly the universal
$R$-matrix. Then we relate them with the chains of coideal
subalgebras of \cite{HS}, and compute the desired
$R$-matrices of quantum doubles of Nichols algebras with finite
root systems. Finally we show some applications of the previous results to relate different PBW bases.

\subsection{PBW bases and Lusztig automorphisms}

Set an element $w=1_\chi s_{i_1}s_{i_2}\cdots s_{i_M}$ of maximal length of
$\cW(\chi)$. Denote
\begin{align}\label{eq:def beta_k}
\beta_k&:= s_{i_1}\cdots s_{i_{k-1}}(\alpha_{i_k}), & 1&\leq k\leq M,
\end{align}
so $\beta_k\neq\beta_l$ if $k\neq l$, and
$\Delta_+^\chi=\{\beta_k|1\leq k\leq M\}$.
Set $q_k:=\chi(\beta_k,\beta_k)$, and $N_k$ the order of
$q_k$, which is possibly infinite. As in \cite[Section 4]{HY-shapov}, set
\begin{align*}
E_{\beta_k}&=T_{i_1}\cdots
T_{i_{k-1}}(E_{i_k})\in\u(\chi)^+_{\beta_k}, &
\overline{E}_{\beta_k}&=T_{i_1}^-\cdots
T_{i_{k-1}}^-(E_{i_k})\in\u(\chi)^+_{\beta_k},
\\ F_{\beta_k}&=T_{i_1}\cdots T_{i_{k-1}}(F_{i_k})\in\u(\chi)^-_{\beta_k}, &
\overline{F}_{\beta_k}&=T_{i_1}^-\cdots
T_{i_{k-1}}^-(F_{i_k})\in\u(\chi)^-_{\beta_k},
\end{align*} for $1\leq k\leq M$.

\begin{theorem} \label{thm: HY PBW bases}
\cite[Theorems 4.5, 4.8, 4.9]{HY-shapov}
The sets
\begin{align*}
& \{ E_{\beta_M}^{a_M}E_{\beta_{M-1}}^{a_{M-1}} \cdots
E_{\beta_1}^{a_1}\, | 0\leq a_k < N_k, \, 1\leq k\leq
M\},
\\ & \{ \overline{E}_{\beta_M}^{a_M}\overline{E}_{\beta_{M-1}}^{a_{M-1}} \cdots \overline{E}_{\beta_1}^{a_1}\, |
0\leq a_k < N_k, \, 1\leq k\leq M\},
\end{align*}
are bases of the vector space $\u^+(\chi)$, and the sets
\begin{align*}
& \{ F_{\beta_M}^{a_M}F_{\beta_{M-1}}^{a_{M-1}} \cdots
F_{\beta_1}^{a_1}\, | 0\leq a_k < N_k, \, 1\leq k\leq
M\},
\\ & \{ \overline{F}_{\beta_M}^{a_M}\overline{F}_{\beta_{M-1}}^{a_{M-1}} \cdots \overline{F}_{\beta_1}^{a_1}\, |
0\leq a_k < N_k, \, 1\leq k\leq M\},
\end{align*}
are bases of the vector space $\u^-(\chi)$. Moreover, for each pair
$1\leq k< l \leq M$,
\begin{align*}
E_{\beta_k}E_{\beta_l}- \chi(\beta_k,\beta_l)
E_{\beta_l}E_{\beta_k} &=
\sum c_{a_{k+1},\ldots,a_{l-1}} E_{\beta_{k+1}}^{a_{k+1}} \cdots
E_{\beta_{l-1}}^{a_{l-1}} \in \u^+(\chi),
\\ \overline{E}_{\beta_k}\overline{E}_{\beta_l}- \chi^{-1}(\beta_k,\beta_l)
\overline{E}_{\beta_l}\overline{E}_{\beta_k} &= \sum
\overline{c}_{a_{k+1},\ldots,a_{l-1}} E_{\beta_{k+1}}^{a_{k+1}}
\cdots E_{\beta_{l-1}}^{a_{l-1}} \in \u^+(\chi),
\\ F_{\beta_k}F_{\beta_l}- \chi(\beta_k,\beta_l)
F_{\beta_l}F_{\beta_k}&= \sum d_{a_{k+1},\ldots,a_{l-1}}
F_{\beta_{k+1}}^{a_{k+1}} \cdots F_{\beta_{l-1}}^{a_{l-1}} \in
\u^-(\chi),
\\ \overline{F}_{\beta_k}\overline{F}_{\beta_l}- \chi^{-1}(\beta_k,\beta_l)
\overline{F}_{\beta_l}\overline{F}_{\beta_k} &= \sum
\overline{d}_{a_{k+1},\ldots,a_{l-1}} F_{\beta_{k+1}}^{a_{k+1}}
\cdots F_{\beta_{l-1}}^{a_{l-1}} \in \u^+(\chi),
\end{align*}
for some $c_{a_{k+1},\ldots,a_{l-1}},
\overline{c}_{a_{k+1},\ldots,a_{l-1}}, d_{a_{k+1},\ldots,a_{l-1}},
\overline{d}_{a_{k+1},\ldots,a_{l-1}}\in\ku$.\qed
\end{theorem}

Note that $E_{\beta_k}E_{\beta_l}- \chi(\beta_k,\beta_l)
E_{\beta_l}E_{\beta_k}= \left[E_{\beta_k}, E_{\beta_l} \right]_c$.

\medskip

Now we want to describe the coproduct of the elements of these PBW
generators. First we introduce the following subspaces of
$\u(\chi)$:
\begin{align*}
B^l_+&:= \langle \{ E_{\beta_l}^{a_l}E_{\beta_{l-1}}^{a_{l-1}} \cdots
E_{\beta_1}^{a_1}\, | 0\leq a_k < N_k\} \rangle \subseteq \u^+(\chi),
\\ C^l_+&:= \langle \{ E_{\beta_M}^{a_M}E_{\beta_{M-1}}^{a_{M-1}} \cdots
E_{\beta_1}^{a_1}\, | \exists j>l \, \mbox{ s.t. } \, a_j\neq 0\} \rangle \subseteq \u^+(\chi),
\\ D^l_+&:= \langle \{ E_{\beta_M}^{a_M}E_{\beta_{M-1}}^{a_{M-1}} \cdots
E_{\beta_1}^{a_1}\, | \exists j<l \, \mbox{ s.t. } \, a_j\neq 0\} \rangle \subseteq \u^+(\chi),
\\ B^l_-&:= \langle \{ F_{\beta_l}^{a_l}F_{\beta_{l-1}}^{a_{l-1}} \cdots
F_{\beta_1}^{a_1}\, | 0\leq a_k < N_k\} \rangle \subseteq \u^-(\chi),
\\ C^l_-&:= \langle \{ F_{\beta_M}^{a_M}F_{\beta_{M-1}}^{a_{M-1}} \cdots
F_{\beta_1}^{a_1}\, | \exists j>l \, \mbox{ s.t. } \, a_j\neq 0\} \rangle \subseteq \u^-(\chi),
\\ D^l_-&:= \langle \{ F_{\beta_M}^{a_M}F_{\beta_{M-1}}^{a_{M-1}} \cdots
F_{\beta_1}^{a_1}\, | \exists j<l \, \mbox{ s.t. } \, a_j\neq 0\} \rangle \subseteq \u^-(\chi),
\end{align*}
$1\leq l\leq M$; $\langle S\rangle$ denotes the subspace spanned by a subset $S$ of $\u(\chi)$.

\begin{prop}\label{prop:B_l coideal subalg}
$B^l_+$ (respectively, $B^l_-$) is a right (respectively, left)
coideal subalgebra of $\u^+(\chi)$ (respectively, $\u^-(\chi)$).
\end{prop}
\pf For each $1\leq l\leq M$, set $w_l=1_\chi s_{i_1}s_{i_{M-1}}\cdots
s_{i_l}$, and the corresponding right coideal subalgebra $B(w_l)$ of
$\u^+(\chi)$ (for the braided coproduct $\underline \Delta$) as in
Theorem \ref{Thm:HSadaptado}; then its Hilbert series is
$$ \cH_{B(w_l)} = \prod_{j=1}^l \bq_{N_l}(X^{\beta_l}).$$
By the definition of $B(w_l)$ in \cite{HS} (which involves the
$T_j$'s) it follows that $E_{\beta_j}\in B(w_l)$ for each $1\leq
j\leq k$. Therefore $B_+^l\subseteq B(w_l)$, because $B(w_l)$ is a
subalgebra. But both $N^\theta_0$-graded vector subspaces of
$\u^+(\chi)$ have the same Hilbert series by Theorem \ref{thm: HY
PBW bases}, so $B^l_+=B(w_l)$ is a right coideal subalgebra.

The statement about $B_-^l$ is analogous because $\u^-(\chi)\simeq
\cB(V^*)^{\cop}$. \epf

\begin{coro}\label{coro:coproduct E_beta}
For each $1\leq l\leq M$,
\begin{align*}
\underline{\Delta}(E_{\beta_l}) & \in E_{\beta_l}\ot 1+1\ot E_{\beta_l} + B^{l-1}_+ \ot C^l_+,
\\ \underline{\Delta}(F_{\beta_l}) & \in F_{\beta_l}\ot 1+1\ot F_{\beta_l} +C^l_+\ot B^{l-1}_-.
\end{align*}
\end{coro}
\pf By the previous Proposition and the fact that $\u^+(\chi)$ is a
graded connected Hopf algebra,
$$ \underline{\Delta}(E_{\beta_l})=E_{\beta_l}\ot 1+1\ot E_{\beta_l} + \sum E_{\beta_{l-1}}^{a_{l-1}} \cdots
E_{\beta_1}^{a_1} \ot X_{a_1,\ldots,a_{l-1}}, $$ for some
$X_{a_1,\ldots,a_{l-1}}\in\u^+(\chi)$. Now express these elements in
terms of the PBW basis:
$$ X_{a_1,\ldots,a_{l-1}}= \sum c_{b_m,\ldots,b_1}^{a_{l-1},\ldots,a_1} E_{\beta_M}^{b_M}E_{\beta_{M-1}}^{b_{M-1}} \cdots
E_{\beta_1}^{b_1}.$$ Suppose that
$c_{b_m,\ldots,b_1}^{a_{l-1},\ldots,a_1}\neq 0$. As $\u^+(\chi)$ is
$\N_0^\theta$-graded, $\beta_l=\sum_{b_i\neq
0}b_i\beta_i+\sum_{a_j\neq 0}a_j\beta_j$. As $j$ runs between 1 and
$l-1$, Theorem \ref{Theorem:equivalenciasconvexo} implies that there
exists $i>l$ such that $b_i\neq 0$.

The proof for $F_{\beta_l}$ is analogous.
\epf

More generally, we can describe the coproduct of each PBW generator.
In this case we can only describe the left hand side of the tensor
product.

\begin{prop}\label{prop:coproduct PBW generators}
For each $1\leq l\leq M$, $1\leq a_l<N_l$,
\begin{align*}
\underline{\Delta}(E_{\beta_l}^{a_l}E_{\beta_{l-1}}^{a_{l-1}} \cdots
E_{\beta_1}^{a_1})  \in &
\sum_{p=0}^{a_l} \binom{a_l}{p}_{q_l} E_{\beta_l}^{p} \ot  E_{\beta_l}^{a_l-p} E_{\beta_{l-1}}^{a_{l-1}} \cdots E_{\beta_1}^{a_1}
\\ & + E_{\beta_l}^{a_l} \cdots E_{\beta_1}^{a_1} \ot 1 + (D^l_+\cap B^l_+)\ot \u^+(\chi) ,
\\ \underline{\Delta}( F_{\beta_l}^{a_l}F_{\beta_{l-1}}^{a_{l-1}} \cdots F_{\beta_1}^{a_1})  \in &
\sum_{p=0}^{a_l} \binom{a_l}{p}_{q_l} F_{\beta_l}^{a_l-p} F_{\beta_{l-1}}^{a_{l-1}} \cdots F_{\beta_1}^{a_1} \ot F_{\beta_l}^{p}
\\ & + 1\ot F_{\beta_l}^{a_l} \cdots F_{\beta_1}^{a_1} + \u^-(\chi) \ot (D^l_-\cap B^l_-) ,\end{align*}
\end{prop}
\pf We prove the statement for the $E_{\beta_k}$'s by induction on
$l$; the proof for the $F_{\beta_k}$'s is analogous. The case $l=1$
is trivial, because $E_{\beta_1}=E_{i_1}$ is primitive, so
$$ \underline{\Delta}(E_{\beta_1}^{a_1}) = \sum_{p=0}^{a_1} \binom{a_1}{p}_{q_1} E_{\beta_1}^{p} \ot E_{\beta_1}^{a_1-p}.$$
Assume that it holds for $k<l$. Now we use induction on $a_l$. If $a_l=1$,
$$ \underline{\Delta}( E_{\beta_l} E_{\beta_{l-1}}^{a_{l-1}} \cdots E_{\beta_1}^{a_1}) = \underline{\Delta}(E_{\beta_l})
\underline{\Delta}(E_{\beta_{l-1}}^{a_{l-1}} \cdots E_{\beta_1}^{a_1}). $$ Therefore we use inductive hypothesis,
Corollary \ref{coro:coproduct E_beta} and the fact that $B_{l-1}$ is
a subalgebra to conclude the proof. The inductive step on $a_l$ is
completely analogous, and close to the proof of results involving
the coproduct of hyperletters in \cite{Kh}. \epf

\subsection{Explicit computation of the universal R-matrix}

We will obtain now an explicit formula for the universal $R$-
matrix when the Nichols algebra is finite-dimensional. By
\eqref{eqn:R-matrix double} it is enough to compute bases of
$\u^{\geq0}(\chi)$ and $\u^{\leq0}(\chi)$, which are dual for
$\eta$. Such bases will be those of Theorem \ref{thm: HY PBW
bases}.

The proof is similar to the one of \cite[Proposition 4.2]{A2}, see
also \cite{R-lyndon}.

\begin{remark}\label{rem:relation coproduct-br. coproduct}
Set for each $\alpha=(a_1,\ldots,a_\theta)\in\Z^\theta$
$$
K^{\alpha}:=K_1^{a_1}\cdots K_\theta^{a_\theta}\in \u^{+0}(\chi),
\qquad L^{\alpha}:=L_1^{a_1}\cdots L_\theta^{a_\theta}\in
\u^{-0}(\chi).$$ For each $\Z^\theta$-homogeneous element
$E\in\u(\chi)$ let $|E|\in\Z^\theta$ be its degree. Therefore,
\begin{equation}\label{eqn:coproduct positive part}
\Delta(E)= E_{(1)}K^{|E_{(2)}|} \otimes E_{(2)},
\end{equation}
for each homogeneous $E\in\u^+(\chi)$, where $\underline{\Delta}(E)=E_{(1)}\ot
E_{(2)}$ is the Sweedler notation for the braided comultiplication.
Analogously, for each homogeneous $F\in\u^-(\chi)$,
\begin{equation}\label{eqn:coproduct negative part}
\Delta(F)= F_{(1)} \otimes F_{(2)}L^{|F_{(1)}|}.
\end{equation}
\end{remark}
\medskip

\begin{prop}\label{prop:duality PBW bases}
Let $0\leq a_i, b_i\leq N_i$, for each $1\leq i\leq M$. Then
\begin{equation}\label{eqn:duality PBW bases}
\eta\left( E_{\beta_M}^{a_M}E_{\beta_{M-1}}^{a_{M-1}} \cdots
E_{\beta_1}^{a_1} , F_{\beta_M}^{b_M}F_{\beta_{M-1}}^{b_{M-1}}
\cdots F_{\beta_1}^{b_1} \right) = \prod_{i=1}^M \delta_{a_i,b_i}
(a_i)_{q_i}! \eta_i^{a_i},
\end{equation}
where $\eta_i:=\eta(E_{\beta_i},F_{\beta_i})$ is not zero for all
$i$.
\end{prop}
\pf We will prove \eqref{eqn:duality PBW bases} by induction on
$\sum a_i$, $\sum b_i$; therefore $\eta_i\neq0$ for all $i$ because $\eta$ is a
non-degenerate pairing. It is clear if $\sum a_i=0$. If $\sum a_i=1$, then the PBW generator is just $E_{\beta_j}$ for some $j$. For this case we apply decreasing induction on $j$.
Note that $\eta\left( E_{\beta_j}, F_{\beta_M}^{b_M}F_{\beta_{M-1}}^{b_{M-1}}
\cdots F_{\beta_1}^{b_1} \right) =0$ when $\beta_j\neq \sum_l b_l\beta_l$, by Proposition \ref{prop:non deg pairing HY}. If $\beta_j= \sum_l b_l\beta_l$ and $\beta_j$ is a simple root, the unique possibility is $b_j=1$ and $b_l=0$ for $l\neq j$. If $\beta_j$ is not a simple root, then either $b_j=1$ and $b_l=0$ for $l\neq j$, or there exists $k>j$ such that $b_k>0$ because the order is strongly convex. In the last case,
\begin{align*}
\eta \left( E_{\beta_j}, F_{\beta_k}^{b_k}F_{\beta_{k-1}}^{b_{k-1}}\cdots F_{\beta_1}^{b_1} \right)&
= \eta \left( (E_{\beta_j})_{(1)}K^{|(E_{\beta_j})_{(2)}|}, F_{\beta_k}\right)
\\ &\eta\left( (E_{\beta_j})_{(2)}, F_{\beta_k}^{b_k-1} F_{\beta_{k-1}}^{b_{k-1}}\cdots F_{\beta_1}^{b_1} \right)=0,
\end{align*}
because $\eta\left( (E_{\beta_j})_{(1)}K^{|(E_{\beta_j})_{(2)}|},
 F_{\beta_k}\right)=0$ by Corollary \ref{coro:coproduct E_beta} and the inductive hypothesis.

Assume that $\sum a_i, \sum b_i>0$ and we have proved the formula for sums smaller than these two. Set $k=\max \{i: a_i\neq 0 \}$, $l=\max \{j: b_j\neq 0 \}$,
and suppose that $k\leq l$ (otherwise the proof is analogous). By Proposition \ref{prop:coproduct PBW generators},
\begin{align*}
\eta & \left( E_{\beta_k}^{a_k}E_{\beta_{k-1}}^{a_{k-1}} \cdots
E_{\beta_1}^{a_1} , F_{\beta_l}^{b_l}F_{\beta_{l-1}}^{b_{l-1}}
\cdots F_{\beta_1}^{b_1} \right)
\\=& \eta\left( (E_{\beta_k}^{a_k}E_{\beta_{k-1}}^{a_{k-1}} \cdots
E_{\beta_1}^{a_1})_{(1)}K^{|(E_{\beta_k}^{a_k}E_{\beta_{k-1}}^{a_{k-1}} \cdots
E_{\beta_1}^{a_1})_{(2)}|} , F_{\beta_l} \right)
\\ &\quad \eta\left( (E_{\beta_1}^{a_1}E_{\beta_{M-1}}^{a_{M-1}} \cdots
E_{\beta_k}^{a_k})_{(2)} , F_{\beta_l}^{b_l-1}F_{\beta_{l-1}}^{b_{l-1}}
\cdots   F_{\beta_1}^{b_1} \right)
\\ =& (b_l)_{q_l}\eta_l\delta_{l,k} \eta\left( E_{\beta_k}^{a_k-1} E_{\beta_{k-1}}^{a_{k-1}}  \cdots
 E_{\beta_1}^{a_1} , F_{\beta_l}^{b_l-1}F_{\beta_{l-1}}^{b_{l-1}}
\cdots   F_{\beta_1}^{b_1} \right),
\end{align*}
so the proof follows by inductive hypothesis.
\epf

Now we obtain a formula for the scalars $\eta_i$. The algebras $\u^{\geq0}(\chi)$, $\u^{\leq0}(\chi)$ are canonically $\N_0$-graded; we denote by $d(X)$, $d(Y)$ the degree of the homogeneous elements $X\in\u^{\geq0}(\chi)$, $Y\in\u^{\leq0}(\chi)$. In fact, if $X\in\u^{\geq0}(\chi)_\beta$, $Y\in\u^{\leq0}(\chi)_{-\beta}$, $\beta=\sum_{i=1}^\theta n_i\alpha_i\in\N_0^\theta$, then $d(X)=d(Y)=\sum_{i=1}^\theta n_i$.

\begin{lemma}
$\eta_k=(-1)^{d(E_{\beta_k})}$ for all $1\leq k\leq M$.
\end{lemma}
\pf By induction on $k$, it is easy to prove that
\begin{equation}\label{eqn:commut EF}
E_{\beta_k}F_{\beta_k}-F_{\beta_k}E_{\beta_k}= K^{\beta_k}-L^{\beta_k}.
\end{equation}
On the other hand, by \eqref{eqn:DbleComXY} we have that
\begin{equation}\label{eqn:EbetaFbeta identity}
E_{\beta_k}F_{\beta_k}= \eta\left((E_{\beta_k})_1,(F_{\beta_k})_1\right)\eta\left((E_{\beta_k})_3,\cS((F_{\beta_k})_3)\right)(F_{\beta_k})_2(E_{\beta_k})_2.
\end{equation}
Using \eqref{eqn:coproduct positive part} and the fact that $\u^{\geq0}(\chi)$ is $\N_0^\theta$-graded, we deduce that the unique term in $\trDeltatwo(E_{\beta_k})$ where appears $K^{\beta_k}$ in the middle is $K^{\beta_k}\otimes K^{\beta_k}\otimes E_{\beta_k}$. If we want to compute the coefficient of in \eqref{eqn:EbetaFbeta identity}, it is enough to look at the term $1\otimes1\otimes F_{\beta_k}$ in $\trDeltatwo(F_{\beta_k})$, because the components of different degree are orthogonal for $\eta$. Using the antipode axiom and that $\u^{\leq0}(\chi)$ is graded, we have that $\cS(F_{\beta_k})$ is written as $(-1)^{d(F_{\beta_k})}F_{\beta_k}L^{-\beta_k}$ plus terms of lower degree. Then the coefficient of $K^{\beta_k}$ in the right hand side of \eqref{eqn:EbetaFbeta identity} is $(-1)^{d(F_{\beta_k})}\eta_k$, using again the orthogonality of the components of different degree.
\epf

\bigbreak

We recall a generalization of Proposition \ref{prop:non deg pairing HY}. The main objective is to consider bosonizations
of Nichols algebras by abelian groups, not only free abelian groups, and their quantum doubles. Similar generalizations
can be found in \cite{ARS,RaS}, and also in \cite{B} for finite groups.

Set a bicharacter $\chi:\Z^\theta\times\Z^\theta\to\ku^\times$, and two abelian groups $\Gamma,\Lambda$. Assume that there exists elements $g_i\in\Gamma$, $\gamma_j\in\widehat{\Gamma}$ such that $\gamma_j(g_i)=q_{ij}$, and elements $h_i\in\Lambda$, $\lambda_j\in\widehat{\Lambda}$ such that $\lambda_j(h_i)=q_{ji}$. Assume that there exists a bicharacter $\mu:\Gamma\times\Lambda\to\ku^\times$, such that $\mu(g_i,h_j)=q_{ij}$. Note that this happens for example when $\Gamma=\Lambda=\Z^\theta$, as in \cite[Section 4]{H-isom}.

Set $V\in\ydg$ as the vector space with a fixed basis $E_1,\ldots,E_\theta$ such that $E_i\in V^{\gamma_i}_{g_i}$, $W\in\ydl$ to the vector space with a fixed basis $F_1,\ldots,F_\theta$ such that $F_i\in V^{\lambda_i}_{h_i}$. Let $\cB=\cB(V)\#\ku\Gamma$ and $\cB'=(\cB(W)\#\ku\Lambda)^{\cop}$.

\begin{theorem}\label{thm:non deg pairing AAB}
There exists a unique skew-Hopf pairing $\eta:\cB\otimes\cB'\to\ku$ such that for all $1\leq i,j\leq\theta$ and all $g\in\Gamma$, $h\in\Lambda$,
\begin{equation}\label{eqn:def eta2}
\eta(g,h)=\mu(g,h), \qquad \eta(E_i,F_j)=-\delta_{ij}, \qquad \eta(E_i,h)=\eta(g,F_j)=0.
\end{equation}
It satisfies the following condition: for all $E\in\u^+(\chi)$, $F\in\u^-(\chi)$, $g\in\cB$, $h\in\cB'$,
\begin{equation}\label{eqn:decomp eta2}
\eta(Eg,Fh)=\eta(E,F)\mu(g,h).
\end{equation}
The restriction of $\eta$ to $\cB(V)\otimes\cB(W)$ coincides with the one of the pairing in Proposition \ref{prop:non deg pairing HY}.\qed
\end{theorem}

We work with the case $\Lambda=\widehat{\Gamma}$, $\Gamma$ a finite group, $\mu$ the evaluation bicharacter, and $h_i=\gamma_i$, $\lambda_i=g_i$ under the canonical identification of the characters of $\widehat{\Gamma}$ with $\Gamma$. In this case $\eta$ is non-degenerate. Call $\u(\chi)$ to the Hopf algebra corresponding to this skew-Hopf pairing, following Subsection \ref{subsection:skew-hopf pairing}, and denote $\cB=\u^{\geq0}(\chi)$, $\cB'=\u^{\leq0}(\chi)$ by analogy with the previous sections. Two dual bases for $\eta|_{\ku\Gamma\otimes\ku\widehat{\Gamma}}$ are $\{g\}_{g\in\Gamma}$, $\{\delta_g\}_{g\in\Gamma}$, where $\delta_g=|\Gamma|^{-1}\sum_{\gamma\in\widehat{\Gamma}}\gamma(g^{-1})\, \gamma$. Therefore it has an $R$-matrix of the form:
\begin{equation}\label{eqn:Rchi}
\cR_1:= \sum_{g\in\Gamma} \delta_g\otimes g= \frac{1}{|\Gamma|}\sum_{g\in\Gamma, \gamma\in
\widetilde{\Gamma}} \gamma(g^{-1}) \, \gamma \otimes g.
\end{equation}

\begin{theorem}\label{thm:R-matrix}
The universal $R$-matrix of $\u(\chi)$ is given by the formula
\begin{equation}\label{eqn:R-matrix}
\cR = \left( \prod \exp_{q_j}\left( (-1)^{d(F_{\beta_k})} F_{\beta_j}\otimes
E_{\beta_j} \right) \right)\cR_1,
\end{equation}
where the product is written in decreasing order.
\end{theorem}
\pf By Proposition \ref{prop:duality PBW bases} and
Theorem \ref{thm:non deg pairing AAB}, the sets
\begin{align*}
& \{ E_{\beta_M}^{a_M}\cdots E_{\beta_1}^{a_1} g: 0\leq a_i< N_i,
g\in\Gamma \},
\\
&\left\{ \left(\prod_{i=1}^M (a_i)_{q_i}! \eta_i^{a_i}\right)^{-1}
\, F_{\beta_M}^{b_M}\cdots F_{\beta_1}^{b_1} \delta_g: 0\leq b_i< N_i,
g\in\Gamma \right\}
\end{align*}
are bases of $\u^{\geq0}(\chi)$, $\u^{\leq0}(\chi)$, respectively,
which are dual for $\eta$. As in Subsection
\ref{subsection:skew-hopf pairing}, a formula for the $R$-matrix is
given by:
\begin{align*}
\cR &= \sum_{g\in\Gamma} \sum_{0\leq a_i< N_i} \left(\prod_{i=1}^M
(a_i)_{q_i}! \eta_i^{a_i}\right)^{-1} \, F_{\beta_M}^{b_M}\cdots
F_{\beta_1}^{b_1} \delta_g \otimes E_{\beta_M}^{a_M}\cdots
E_{\beta_1}^{a_1} g
\\ &= \left( \prod \left( \sum_{i=0}^{N_j-1} \frac{\eta_j^i}{(i)_{q_j}!} F_{\beta_j}^i\otimes E_{\beta_j}^i
\right) \right) \left(\sum_{g\in\Gamma} \delta_g\otimes g \right),
\end{align*}
which ends the proof. \epf

\bigskip

\subsection{Further computations on convex PBW bases}

We can refine the coproduct expression of each $E_\beta$.
In consequence we can obtain a family of left coideal subalgebras,
induced by products of the same PBW generators. For each $1\leq l\leq M$, let
\begin{align*}
\mathbf{B}^l_+&:= \langle \{ E_{\beta_M}^{a_M}E_{\beta_{M-1}}^{a_{M-1}} \cdots
E_{\beta_l}^{a_l}\, | 0\leq a_k < N_k\} \rangle \subseteq \u^+(\chi),
\\ \mathbf{B}^l_-&:= \langle \{ F_{\beta_M}^{a_M}F_{\beta_{M-1}}^{a_{M-1}} \cdots
F_{\beta_l}^{a_l}\, | 0\leq a_k < N_k\} \rangle \subseteq \u^-(\chi).
\end{align*}

\begin{lemma}\label{lemma:coproduct E_beta refined}
For each $1\leq l\leq M$,
\begin{align*}
\underline{\Delta}(E_{\beta_l}) & \in E_{\beta_l}\ot 1+1\ot E_{\beta_l} + B^{l-1}_+ \ot \mathbf{B}^{l-1}_+,
\\ \underline{\Delta}(F_{\beta_l}) & \in F_{\beta_l}\ot 1+1\ot F_{\beta_l} + \mathbf{B}^{l-1}_-\ot B^{l-1}_-.
\end{align*}
\end{lemma}
\pf
Write both sides of $\underline{\Delta}(E_{\beta_l})$ as linear combinations of the elements of the PBW basis, and take a term
$$ E_{\beta_{l-1}}^{a_{l-1}} \cdots E_{\beta_1}^{a_1} \otimes E_{\beta_M}^{b_M}E_{\beta_{M-1}}^{b_{M-1}} \cdots E_{\beta_k}^{b_k} $$
appearing with non-zero coefficient $c$, where $k$ is such that $b_k\neq 0$. Using the orthogonality of the elements of the PBW basis,
\begin{align*}
 0\neq c\, \eta&(E_{\beta_{l-1}}^{a_{l-1}} \cdots E_{\beta_1}^{a_1}K^{|E_{\beta_M}^{b_M}E_{\beta_{M-1}}^{b_{M-1}} \cdots E_{\beta_k}^{b_k}|}, F_{\beta_{l-1}}^{a_{l-1}} \cdots F_{\beta_1}^{a_1})
\\ &\eta(  E_{\beta_M}^{b_M}E_{\beta_{M-1}}^{b_{M-1}} \cdots E_{\beta_k}^{b_k} ,  F_{\beta_M}^{b_M}F_{\beta_{M-1}}^{b_{M-1}} \cdots F_{\beta_k}^{b_k} )
\\ = \eta&(E_{\beta_l}, F_{\beta_{l-1}}^{a_{l-1}} \cdots F_{\beta_1}^{a_1} F_{\beta_M}^{b_M}F_{\beta_{M-1}}^{b_{M-1}} \cdots F_{\beta_k}^{b_k}).
\end{align*}
Suppose that $k<l$. Using last part of Theorem \ref{thm: HY PBW bases} repeatedly we see that
$$ z:=F_{\beta_{l-1}}^{a_{l-1}} \cdots F_{\beta_1}^{a_1} F_{\beta_M}^{b_M}F_{\beta_{M-1}}^{b_{M-1}} \cdots F_{\beta_k}^{b_k} \in D^l_-,$$
so $\eta(E_{\beta_l},z)=0$, a contradiction. Then $k\geq l$, and we end the proof.
\epf

\begin{prop}\label{prop:B_l coideal subalg II}
$\mathbf{B}^l_+$ (respectively, $\mathbf{B}^l_-$) is a left (respectively, right)
coideal subalgebra of $\u^+(\chi)$ (respectively, $\u^-(\chi)$).
\end{prop}
\pf
It is a consequence of Lemma \ref{lemma:coproduct E_beta refined} and last part of Theorem \ref{thm: HY PBW bases}.
\epf

For the last part of this section we prove a result generalizing \cite[Theorem 22]{R-lyndon}. It establishes the uniqueness
(up to scalars) of a PBW basis determining a filtration of coideal subalgebras, and it is useful to compare PBW bases
coming from Lusztig isomorphisms as in the previous results, and PBW bases from combinatorics as \cite{Kh}. Note that the
first kind of PBW bases gives right and left coideal subalgebras, while some examples of the second family give left coideal
subalgebras, see \cite[Section 3.3]{A2}.

\begin{theorem}\label{thm: comparing bases PBW}
 Let $(\En_\beta)_{\beta\in\Delta_+^\chi}$ be non-zero elements of $\u^+(\chi)$, such that $\En_\beta\in\u^+(\chi)_{\beta}$, and there exists an order $\beta_M>\ldots >\beta_1$ on the roots such that, for each $1\leq k\leq M$, the elements $\En_{\beta_M}^{a_M}\cdots \En_{\beta_k}^{a_k}$, $0\leq a_j< N_{\beta_k}$, determine a basis of a subspace $\mathbf{Y}_k$, which is a left coideal subalgebra of $\u^+(\chi)$. Then the order on the roots is convex.

Moreover, if $(E_{\beta})_{\beta\in\Delta_+^\chi}$ denote PBW generators for the corresponding expression of the element of maximal length of $\cW$, then there exists non-zero scalars $c_{\beta}$ such that $\En_\beta=c_\beta\, E_\beta$.
\end{theorem}
\pf
The convexity on the order follows from the fact that the chain of coideal subalgebras $\mathbf{Y}_M\subsetneq \cdots\subsetneq\mathbf{Y}_1=\cB(V)$ coincides with $\mathbf{B}_+^M\subsetneq \cdots\subsetneq\mathbf{B}_+^1=\cB(V)$. The proof of this fact is exactly as in \cite[Theorem 3.16]{A2}. That is, $\mathbf{Y}_k = \mathbf{B}_+^k$ for all $1\leq k\leq M$..

For the second statement, write $\En_{\beta_k}= \sum c(a_1,\cdots,a_M)\, E_{\beta_M}^{a_M} \cdots E_{\beta_1}^{a_1}$. If $c(a_1,\cdots,a_M)\neq 0$, then $\beta_k=\sum_j a_j\beta_j$, so $a_k=1$, $a_j=0$ for all $j\neq k$, or there exists $j<k$ such that $a_j\neq0$. The second case is not possible because $\En_{\beta_k}\in\mathbf{Y}_k = \mathbf{B}_+^k$. Therefore, $\En_{\beta_k}=c_{\beta_k}\, E_{\beta_k}$ for some $c_{\beta_k}\in\ku^\times$.
\epf

\begin{example}
Let $(q_{ij})_{1\leq i,j\leq 2}$ be a matrix whose generalized Dynkin diagram is $\xymatrix{\circ^{\zeta} \ar@{-}[r]^{\zeta^2} & \circ^{-1}}$, $\zeta$ a root of unity of order $5$, and $\chi$ the associated bicharacter on $\Z^2$. The element of maximal length on its Weyl groupoid has a reduced expression $w_0=\id^\chi s_1s_2s_1s_2s_1s_2s_1s_2$. Then
\begin{align*}
\alpha_1&<3\alpha_1+\alpha_2<2\alpha_1+\alpha_2<5\alpha_1+3\alpha_2
\\ &<3\alpha_1+2\alpha_2<4\alpha_1+3\alpha_2<\alpha_1+\alpha_2<\alpha_2
\end{align*}
is the corresponding order on the roots. We obtain a PBW basis with generators $E_\beta$, $\beta\in\Delta_+^\chi$, using the Lusztig isomorphisms. Let $\Gamma$ be a finite abelian group, $g_1,g_2\in\Gamma$, $\gamma_1,\gamma_2\in\widehat\Gamma$ such that $\gamma_j(g_i)=q_{ij}$, so $\cB(V)$ can be viewed as a braided Hopf algebra on the category of Yetter-Drinfeld modules of $\ku\Gamma$. We define $\cR_1$ as in \eqref{eqn:Rchi}.
By Theorem \ref{thm:R-matrix},
\begin{align*}
\cR &= 
\left( \sum_{k=0}^4 \frac{-1}{(k)_\zeta !}F_1\otimes E_1 \right)
\left( 1\ot 1 - F_{3\alpha_1+\alpha_2}\otimes E_{3\alpha_1+\alpha_2}  \right) \\
&\qquad \left( \sum_{k=0}^9 \frac{-1}{(k)_{-\zeta^3} !}F_{2\alpha_1+\alpha_2}\otimes E_{2\alpha_1+\alpha_2} \right)
\left( 1\ot 1 + F_{5\alpha_1+3\alpha_2}\otimes E_{5\alpha_1+3\alpha_2}  \right) \\
&\qquad \left( \sum_{k=0}^4 \frac{-1}{(k)_\zeta !}F_{3\alpha_1+2\alpha_2}\otimes E_{3\alpha_1+2\alpha_2} \right)
\left( 1\ot 1 - F_{4\alpha_1+3\alpha_2}\otimes E_{4\alpha_1+3\alpha_2}  \right) \\
&\qquad \left( \sum_{k=0}^9 \frac{-1}{(k)_{-\zeta^3} !}F_{\alpha_1+\alpha_2}\otimes E_{\alpha_1+\alpha_2} \right)
\left( 1\ot 1 - F_2\otimes E_2  \right)  \cR_1.
\end{align*}

We can obtain also a PBW basis of hyperletters $\En_\beta=[\ell_\beta]_c$, $\beta\in\Delta_+^\chi$, associated to Lyndon words $\ell_\beta$ as in \cite{Kh}. We compute easily the corresponding Lyndon words using \cite[Corollary 3.17]{A2}:
\begin{align*}
\ell_{\alpha_1}&=x_1, & \ell_{3\alpha_1+\alpha_2}&=x_1^3x_2, & \ell_{4\alpha_1+3\alpha_2}&=x_1^2x_2x_1x_2x_2x_1x_2, \\
\ell_{\alpha_2}&=x_2, & \ell_{2\alpha_1+\alpha_2}&=x_1^2x_2, & \ell_{5\alpha_1+3\alpha_2}&=x_1^2x_2x_1^2x_2x_1x_2,  \\
\ell_{\alpha_1+\alpha_2}&=x_1x_2, & \ell_{3\alpha_1+2\alpha_2}&=x_1^2x_2x_1x_2.
\end{align*}
We compute using the Shirshov decomposition, see \cite{A2, Kh} and the references there in,
\begin{align*}
\En_{\alpha_1}&=x_1, & \En_{3\alpha_1+\alpha_2}&=(\ad_cx_1)^3x_2,  \\
\En_{\alpha_2}&=x_2, & \En_{3\alpha_1+2\alpha_2}&=[\En_{2\alpha_1+\alpha_2}, \En_{\alpha_1+\alpha_2}]_c, \\
\En_{\alpha_1+\alpha_2}&=(\ad_cx_1)x_2, & \En_{4\alpha_1+3\alpha_2}&=[\En_{3\alpha_1+2\alpha_2}, \En_{\alpha_1+\alpha_2}]_c,\\
\En_{2\alpha_1+\alpha_2}&=(\ad_cx_1)^2x_2, & \En_{5\alpha_1+3\alpha_2}&=[\En_{2\alpha_1+\alpha_2}, \En_{3\alpha_1+2\alpha_2}]_c.
\end{align*}
By the previous theorem, there exists $c_{\beta}\in\ku^\times$ such that $\En_\beta=c_\beta\, E_\beta$. It can be computed as the inverse of the coefficient
of $\ell_\beta$ in $E_\beta$, because $\ell_\beta$ appears with coefficient $1$ in $\En_\beta$.
\end{example}

\end{document}